\ifx\fltex\undefined
\documentclass[12pt,a4paper,english]{article}
\usepackage{amsmath,amssymb}
\usepackage{babel}
\usepackage{eucal}
\fi 

\DeclareMathVersion{scriptstyle}

\usepackage{isolatin1,theorem}
\theoremstyle{change}
{\theorembodyfont{\slshape} 
  \newtheorem{thm}{Theorem.}[section]
  \newtheorem{lemma}[thm]{Lemma.}
  \newtheorem{prop}[thm]{Proposition.}
  \newtheorem{cor}[thm]{Corollary.}
}
{\theorembodyfont{\rmfamily} 
  \newtheorem{definition}[thm]{Definition.}
  \newtheorem{remark}[thm]{Remark.}

}

\numberwithin{equation}{section} 

\newcommand{\proof}[1][Proof. ]{{\it#1}}

\def\endproof{{\nobreak\qquad{\scriptstyle \blacksquare}}}

\def\C{{\mathbb C}}
\def\B{{\mathbb B}}

\def\F{{\mathbb F}}

\def\N{{\mathbb N}}
\def\R{{\mathbb R}}
\def\T{{\mathbb T}}

\def\CB{{\mathcal B}}

\def\CD{{\mathcal D}}

\def\CH{{\mathcal H}}

\def\CJ{{\mathcal J}}

\def\CM{{\mathcal M}}
\def\CN{{\mathcal N}}

\def\CR{{\mathcal R}}
\def\CS{{\mathcal S}}

\def\CU{{\mathcal U}}
\def\CV{{\mathcal V}}

\def\unit{{\bf 1}}
\def\i{{\rm i}}
\def\e{{\rm e}}
\def\d{{\rm d}}
\def\eps{\varepsilon}
\def\<{{\langle}}
\def\>{{\rangle}}

\def\CMD{\CM^\Delta}

\def\im{{\text{\rm Im}}}
\def\re{{\text{\rm Re}}}
\def\supp{{\text{\rm supp}}}

\def\Prob{{\text{\rm Prob}}}

\def\Lp{L^p(\CM, \tau)}

\def\cc{^*}

\begin{document}

\title{{\Large{\bf Brown Measures of Unbounded Operators \\ Affiliated with a Finite von Neumann Algebra}}}
\author{Uffe~Haagerup and Hanne~Schultz \footnote{Supported by the Danish~National~Research~Foundation.}}      
\maketitle

\begin{abstract}
\noindent In this paper we generalize Brown's spectral distribution measure to
  a large class of unbounded operators affiliated with a finite von
  Neumann algebra. Moreover, we compute the Brown measure of all
  unbounded $R$--diagonal operators in this class. As a particular
  case, we determine the Brown measure $z=xy^{-1}$, where $(x,y)$ is a
  circular system in the sense of Voiculescu, and we prove that for
  all $n\in\N$, $z^n\in\Lp$ if and only if $0<p<\frac{2}{n+1}$. 
\end{abstract}

\section{Introduction}

Let $\CM$ be a finite von Neumann algebra with a faithful, normal,
tracial state $\tau$, and let
\[
\Delta(T)=\exp \Big(\int_0^\infty \log t\,\d\mu_{|T|}(t)\Big)
\]
denote the corresponding Fuglede--Kadison determinant. L.~G.~Brown
proved in \cite{Bro} that for every $T\in\CM$, there exists a unique,
compactly supported measure $\mu_T\in {\rm Prob}(\C)$ with the
property that
\[
\log\Delta(T-\lambda\unit)=\int_\C \log|z-\lambda|\,\d\mu_T(z), \qquad
\lambda\in\C.
\]
This measure is called Brown's spectral distribution measure (or just
the Brown measure) of $T$. It was computed in a number of special
cases in \cite{HL}, \cite{BL}, \cite{DH3}, and \cite{AH}. In particular, it was proven in
\cite[Theorem~4.5]{HL} that if $T\in\CM$ is $R$--diagonal in the sense
of Nica and Speicher \cite{NS}, then $\mu_T$ can be determined from
the $S$--transform of the distribution $\mu_{|T|^2}$. For simplicity,
assume that $T\in\CM$ is an $R$--diagonal element which is not
proportional to a unitary and for which ${\rm ker}(T)=0$. Then $\mu_T$
is the unique probability measure on $\C$ which is invariant under the
rotations $z\mapsto \gamma z$, $\gamma\in\T$, and which satisfies
\[
\mu_T\Big(B\big(0,\CS_{\mu_{|T|^2}}(t-1)^{-\frac12}\big)\Big) = t,
\qquad 0<t<1.
\]
In this paper we extend the Brown measure to all operators in the set
$\CMD$ of closed, densely defined operators $T$ affiliated with $\CM$
satisfying
\[
\int_0^\infty \log^+ t\,\d\mu_{|T|}(t)<\infty.
\]
Moreover, we extend \cite[Theorem~4.5]{HL} to all $R$--diagonal
operators in $\CMD$. Finally, we will study a particular example of an
unbounded $R$--diagonal element, namely the operator $z=xy^{-1}$,
where $(x,y)$ is a circular system in the sense of Voiculescu.

The material in this paper is organized as follows: In section~2 we
introduce the class $\CMD$ and generalize the Brown measure to all
$T\in\CMD$ by proving, that for such $T$, there is a unique $\mu_T\in
{\rm Prob}(\C)$ satisfying
\[
\int_\C\log^+|z|\,\d\mu_T(z) < \infty
\]
and
\[
\log\Delta(T-\lambda\unit)=\int_\C\log|z-\lambda|\,\d\mu_T(z), \qquad
\lambda\in\C.
\]
Moreover, we extend Weil's inequality
\[
\int_\C |z|^p\,\d\mu_T(z)\leq \|T\|_p^p
\]
to all $T\in\Lp$. The main results in section~2 are stated in
the appendix of Brown's paper \cite{Bro} without proofs or with very sketchy
proofs. Since the results of the remaining sections of this paper and
of our forthcoming paper \cite{HS2} rely heavily on these statements,
we have decided to include complete proofs. We will follow a different
route than the one outlined in \cite{Bro}. For instance, we do not use the
functions $\Lambda_t(T)$ and $s_T(t)$ from \cite[section~1]{Bro}.

In section~3 we introduce unbounded $R$--diagonal operators and we
prove the following generalization of \cite[section~3]{HL}: The powers
$(S^n)_{n=1}^\infty$ of an $R$--diagonal operator are $R$--diagonal,
and the sum $S+T$ and the product $ST$ of $\ast$--free $R$--diagonal
operators are again $R$--diagonal. Moreover,
\begin{eqnarray*}
  \mu_{|S^n|^2}&=&\mu_{|S|^2}^{\boxtimes n},\\
  \tilde{\mu}_{|S+T|}&=& \tilde\mu_{|S|}\boxplus \tilde\mu_{|T|},\\
  \mu_{|ST|^2}&=& \mu_{|S|^2}\boxtimes \mu_{|T|^2},
\end{eqnarray*}
where $\tilde\mu = \frac12 (\mu + \check{\mu})$ denotes the
symmetrization of a measure $\mu\in {\rm Prob}(\R)$, and $\boxplus$
($\boxtimes$, resp.) denotes the additive (multiplicative, resp.) free
convolution of measures (cf. \cite{BV}). These results are applied in
section~4 to determine the Brown measure of $R$--diagonal operators in
$\CMD$.

In section~5 we consider the operator $z=xy^{-1}$, where $(x,y)$ is a
circular system in the sense of Voiculescu, and we prove that the
Brown measure of $z$ is given by
\[
\d\mu_z(s)= \frac{1}{\pi(1+|s|^2)}\,\d\re s\,\d \im s.
\]
Moreover, we show that for all $n\in\N$, $z^n, z^{-n}\in\Lp$ iff
$0<p<\frac{2}{n+1}$, and when this holds,
\[
\|z^n\|_p^p = \|z^{-n}\|_p^p = \frac{(n+1)\sin\Big(\frac{\pi
    p}{2}\Big)}{\sin\Big(\frac{(n+1)\pi p}{2}\Big)},
\]
and
\[
\|(z^n-\lambda\unit)^{-1}\|_p\leq \|z^{-n}\|_p, \qquad \lambda \in \C.
\]
The last two formulas play a key role in our forthcoming paper
\cite{HS2} on invariant subspaces for operators in a general
II$_1$--factor. 
\section{The Brown measure of certain unbounded operators.}

In \cite[Appendix]{Bro} Brown described in outline how to define a Brown
measure for certain {\it undbounded} operators affiliated with $\CM$, where
$\CM$ is a von Neumann algebra equipped with a faithful, normal, semifinite trace.

In this section we give a more detailed exposition on the
subject in the case where $\CM$ is a finite von Neumann algebra with faithful, tracial
state $\tau$. To be more explicit, we show how one can extend the
definition of the Brown measure to a class $\CMD$ of closed, densely
defined operators affiliated with $\CM$. We also prove that many of the
properties of the Brown measure for bounded operators carry over to the
unbounded case. 

We let $\tilde\CM$ denote the set of closed, densely defined
operators affiliated with $\CM$. Recall that every
operator $T\in\tilde\CM$ has a polar decomposition
\begin{equation}\label{eq1-1}
  T= U|T| = U\int_0^\infty t\,\d E_{|T|}(t),
\end{equation}

where $U\in\CM$ is a unitary, and the spectral measure $E_{|T|}$ takes
values in $\CM$. In particular, for $T\in\tilde\CM$ we may define $\mu_{|T|}\in\Prob(\R)$ by
\begin{equation}\label{eq1-2}
  \mu_{|T|}(B)= \tau(E_{|T|}(B)), \qquad (B\in\B).
\end{equation}

\begin{definition} We denote by $\CM^\Delta$ the set of operators
  $T\in\tilde\CM$ fulfilling the condition
  \begin{equation}\label{eq1-3}
    \tau(\log^+|T|)= \int_0^\infty \log^+(t)\,\d\mu_{|T|}(t)<\infty.
  \end{equation}
  For $T\in\CM^\Delta$, the integral
\[
\int_0^\infty \log t\,\d\mu_{|T|}(t) \,\in \,\R\cup\{-\infty\}
\]
is well--defined, and we define the {\it Fuglede--Kadison determinant}
of $T$, $\Delta(T)\in [0,\infty)$, by
\begin{equation}
\Delta(T)= \exp\Big(\int_0^\infty \log t\,\d\mu_{|T|}(t)\Big).
\end{equation}
\end{definition}

Note that for $T\in\CM$, $\Delta(T)$ is the usual Fuglede--Kadison
determinant of $T$.

\vspace{.2cm}

\begin{remark}
  If $T\in \Lp$ for some $p\in (0,\infty)$, then
    \[
    \int_0^\infty t^p\,\d \mu_{|T|}(t)<\infty,
    \]

    implying that
    \[
    \int_1^\infty \log t\,\d \mu_{|T|}(t) = \frac 1p \int_1^\infty \log (t^p)\,\d
    \mu_{|T|}(t) \leq \frac 1p \int_1^\infty t^p\,\d \mu_{|T|}(t) <\infty,
    \]

    and hence $T\in\CM ^\Delta$.
\end{remark}

\vspace{.2cm}

\begin{lemma}\label{lemma2.3} If $T\in\CM^\Delta$ and $\Delta(T)>0$, then $T$ is
  invertible in $\tilde\CM$, $T^{-1}\in\CM^\Delta$, and
  $\Delta(T^{-1})= \frac{1}{\Delta(T)}$.
\end{lemma}

\proof If $T\in\CM^\Delta$ and $\Delta(T)>0$, then
\[
\int_0^1 |\log t|\,\d\mu_{|T|}(t) <\infty.
\]
Hence, $\tau(E_{|T|}(\{0\}))=\mu_{|T|}(\{0\})=0$, so that $\ker(T)=\{0\}$. Since $\CM$ is
finite, also $\ker(T\cc)=\{0\}$, which implies that $T$ has a closed,
densely defined inverse $T^{-1}\in\tilde\CM$. Take a unitary $U\in\CM$
such that $T=U|T|$. Then
\[
|T^{-1}| = U|T|^{-1}U\cc.
\]
Hence, $\mu_{|T^{-1}|}= \mu_{|T|^{-1}}$. Since $\mu_{|T|^{-1}}$ is the
push--forward measure of $\mu_{|T|}$ via the map $t\mapsto \frac 1t$,
we now have that
\begin{eqnarray*}
\int_1^\infty \log t\,\d\mu_{|T^{-1}|}(t) &=& \int_1^\infty
\log t\,\d\mu_{|T|^{-1}}(t)\\
&=& \int_0^1\log\Big(\frac 1t\Big) \,\d\mu_{|T|}(t)\\
&=&-\int_0^1\log t\,\d\mu_{|T|}(t)\\
&<& \infty.
\end{eqnarray*}
Hence, $T^{-1}\in\CM^\Delta$ and
\[
\log \Delta(T^{-1}) = \int_0^\infty \log\Big(\frac
1t\Big)\,\d\mu_{|T|}(t) = -\log \Delta(T),
\]
i.e. $\Delta(T^{-1})= \frac{1}{\Delta(T)}$. $\endproof$



\vspace{.2cm}

\begin{lemma} Let $T\in\tilde\CM$. Then the following are equivalent:
  \begin{itemize}
    \item[(a)] $T\in\CM^\Delta$, i.e. $\int_0^\infty
    \log^+(t)\,\d\mu_{|T|}(t)<\infty$.
    \item[(b)] $T=A B^{-1}$ for some $A, B\in\CM$ with $\Delta(B)>0$.
    \item[(c)] $T= C^{-1} D$ for some $C, D\in\CM$ with $\Delta(C)>0$.  
  \end{itemize}
  Moreover, if $T\in\CM^\Delta$ and $T=A B^{-1}=C^{-1}D$ for some $A,
  B,C,D\in\CM$ with $\Delta(B),\Delta(C)>0$, then
  \begin{equation}\label{Delta abcd}
    \Delta(T)=\frac{\Delta(A)}{\Delta(B)} =
    \frac{\Delta(D)}{\Delta(C)}.
  \end{equation}
\end{lemma}


\proof If $T\in\CMD$, then $T=U|T|$ for some unitary $U\in \CM$, and $T=AB^{-1}$, where
\begin{equation}\label{a}
  A=U |T|(|T|^2+\unit)^{-\frac 12}\in \CM
\end{equation}
and
\begin{equation}\label{b}
  B = (|T|^2+\unit)^{-\frac 12} \in \CM.
\end{equation}

 Since $\frac 12 \log(t^2+1)\leq \log(2t)$ when $t\geq
1$, we get that
\begin{eqnarray}\label{eq1-11}
  \log \Delta(B) &=& -\frac 12 \int_0^\infty
                 \log(t^2+1)\,\d\mu_{|T|}(t)\nonumber \\
                 &\geq & -\frac 12 \int_{[0,1[}\log 2\,\d\mu_{|T|}(t) -
                 \int_{[1,\infty[}\log(2t)\d\mu_{|T|}(t) \nonumber \\
                 & >& -\infty,
\end{eqnarray}

that is, $\Delta(B) >0$.

Also, $T= U|T|U\cc U$, and with
\begin{equation}
  S= U|T|U\cc,
\end{equation}
\begin{equation}
C = (S^2+\unit)^{-\frac 12}\in \CM,
\end{equation}
and
\begin{equation}
D = S (S^2+\unit)^{-\frac 12}\in \CM,
\end{equation}
we have that $T= C^{-1}DU$. Moreover,
\begin{eqnarray*}
\log \Delta(C) &=& -\frac 12 \int_0^\infty \log(t^2+1)\,\d\mu_{S}(t)\\
&=&  -\frac 12 \int_0^\infty \log(t^2+1)\,\d\mu_{|T|}(t)\\
&>& -\infty,
\end{eqnarray*}
i.e. $\Delta(C)>0$.

Now we have shown that (a) implies (b) and (c). On the other hand, if $T =AB^{-1}$ for some $A, B\in\CM$ with
$\Delta(B)>0$, then we may assume that $B\geq 0$. Then
\begin{eqnarray*}
  \tau(\log^+|T|) &\leq & \tau(\log(\unit + |T|^2))\\
  &=& \tau(\log(\unit + B^{-1}A\cc A B^{-1})).
\end{eqnarray*}
Since $B^{-1}A\cc A B^{-1}\leq \|A\|^2B^{-2}$, and since $t\mapsto
\log(1+t)$ is operator monotone on $[0,\infty)$, we get that
\begin{eqnarray*}
   \tau(\log^+|T|) &\leq & \tau(\log(\unit + \|A\|^2B^{-2}))\\
   & \leq & \tau(\log((1+\|A\|^2)(\unit + B^{-2})))\\
   &=& \log(1+\|A\|^2)+ \tau(\log(\unit + B^{-2})).
\end{eqnarray*}
Since $B$ is bounded and $\Delta(B)>0$,
\begin{eqnarray*}
  \tau(\log(\unit + B^{-2})) &=& \tau(\log(B^2 +\unit))-2\tau(\log
  B)\\
  &\leq & \log(\|B\|^2 +1) - 2\Delta(B)\\
  &<& \infty.
\end{eqnarray*}
This shows that $T\in \CMD$, i.e. (b) implies (a). It follows that if $T= C^{-1} D$ for some $C,
D\in\CM$ with $\Delta(C)>0$, then $T\cc \in \CMD$. Take a unitary $U\in\CM$ such
that $T=U|T|$. Then $|T\cc| = U|T|U\cc$, implying that
$\mu_{|T\cc|} = \mu_{|T|}$. Hence $T$ belongs to $\CMD$ as well, and (c)
implies (a).

Now, let $T\in\CMD$. Then $T=AB^{-1}=C^{-1}D$ for some  $A,B,C,D\in\CM$ with
$\Delta(B),\Delta(C)>0$. Moreover, for all such choices of $A,B,C$ and
$D$,
\[
CA = C(AB^{-1})B = C(C^{-1}D)B = DB.
\]
Since $\Delta$ is multiplicative on $\CM$ (cf. \cite{FuKa}), it follows
that
\[
\Delta(C)\Delta(A)= \Delta(CA)=\Delta(DB)=\Delta(D)\Delta(B).
\]
Hence,
\begin{equation}\label{broeker}
\frac{\Delta(A)}{\Delta(B)}=\frac{\Delta(D)}{\Delta(C)}.
\end{equation}
In particular, with $A$ and $B$ as in \eqref{a} and \eqref{b},
respectively, we have that $\Delta(B)>0$, $T=AB^{-1}$, and 
\[
\log\Delta(A)= \int_0^\infty
\log\Bigg(\frac{t}{\sqrt{t^2+1}}\Bigg)\,\d\mu_{|T|}(t),
\]
and
\[
\log\Delta(B) = \int_0^\infty\log\Bigg(\frac{1}{\sqrt{t^2+1}}\Bigg)\,\d\mu_{|T|}(t), 
\]
so that
\[
\log\Delta(T)= \log\Delta(A)-\log\Delta(B).
\]
Then by \eqref{broeker}, for all choices of $C,D\in\CM$ with
$\Delta(C)>0$ and $T=C^{-1}D$,
\[
\frac{\Delta(D)}{\Delta(C)}=\frac{\Delta(A)}{\Delta(B)}= \Delta(T).
\]
Then finally, by \eqref{broeker}, for all choices of $A,B\in\CM$ with
$\Delta(B)>0$ and $T=AB^{-1}$, we also have that
\[
\frac{\Delta(A)}{\Delta(B)}= \Delta(T). \quad \endproof
\]

\vspace{.2cm}

\begin{prop}\label{prop2.5} If $S,T\in\CMD$, then $ST\in\CMD$, and
\begin{equation}\label{Deltamult}
  \Delta(ST)=\Delta(S)\Delta(T).
\end{equation}
\end{prop}

\proof Let $S,T\in\CMD$. Take $A,B,C,D\in\CM$ with
$\Delta(B),\Delta(C)>0$, such that $T=AB^{-1}$ and $S=C^{-1}D$. Then
\[
ST= C^{-1} DAB^{-1},
\]
where $DAB^{-1}\in \CMD$. Hence there exist $E, F\in\CM$ with
$\Delta(E)>0$ such that  $DAB^{-1} = E^{-1} F$. It follows
that
\begin{equation}\label{e,f}
ST =  C^{-1}E^{-1} F= (EC)^{-1} F,
\end{equation}
where $EC, F\in \CM$, and $\Delta(EC)=\Delta(E)\Delta(C)>0$. That is, $ST$
belongs to $\CMD$.

To prove \eqref{Deltamult}, we let $A,B,C,D,E,F$ be as above. Applying
\eqref{Delta abcd} to $ST=(EC)^{-1}F$, $S=C^{-1}D$ and $T=AB^{-1}$, we
get that
\[
\Delta(ST)=\frac{\Delta(F)}{\Delta(EC)} =
\frac{\Delta(F)}{\Delta(E)\Delta(C)} =
\frac{\Delta(DA)}{\Delta(B)}\frac{1}{\Delta(C)} =
\frac{\Delta(A)}{\Delta(B)}\frac{\Delta(D)}{\Delta(C)} =
\Delta(S)\Delta(T). \quad \endproof
\]

\vspace{.2cm}

\begin{prop}\label{prop2.5.5} $\CMD$ is a subspace of $\tilde\CM$. In
  particular, for $T\in\CMD$ and $\lambda\in\C$,
  $T-\lambda\unit\in\CMD$.
\end{prop}

\proof Clearly, if $T\in\CMD$ and $\alpha\in\C$, then $\alpha
T\in\CMD$. If $S,T\in\CM$, choose $A,B,C,D\in\CM$ with $\Delta(B)>0$,
$\Delta(C)>0$ and such that
\[
S=C^{-1}D, \qquad T= AB^{-1}.
\]
Then
\[
S+T= C^{-1}(DB + CA)B^{-1},
\]
where $DB+CA\in\CM$ and $B^{-1},C^{-1}\in\CMD$
(cf. Lemma~\ref{lemma2.3}). Then, by Proposition~\ref{prop2.5},
$S+T\in\CMD$. $\endproof$
\vspace{1cm}

In the following we consider a fixed operator $T\in \CMD$. Then we define
$f: \C\rightarrow [-\infty,\infty)$ by
\begin{equation}\label{eq1-4}
  f(\lambda) = L(T-\lambda\unit):= \log\Delta(T-\lambda\unit), \qquad (\lambda\in\C).
\end{equation}

\vspace{.2cm}

The next thing we want to prove is:

\begin{thm}\label{mainthm} $f$ given by \eqref{eq1-4} is subharmonic in $\C$, and
  \begin{equation}\label{eq1-5}
    \d\mu_T= \frac{1}{2\pi}\nabla^2 f \,\d\lambda
  \end{equation}
  (taken in the distribution sense) defines a probability measure on $(\C,
  \B_2)$. $\mu_T$ is the unique probability measure on $(\C,\B_2)$
  satisfying
\begin{itemize}
    \item[(i)]
        \[
        \int_\C\log^+|z|\,\d\mu_T(z)<\infty,
        \]
    \item[(ii)] 
       \begin{equation}\label{eq1-7}
        \forall \lambda\in\C:\quad L(T-\lambda\unit)= \int_\C \log|\lambda-z|\,\d \mu_T(z).
      \end{equation}
  \end{itemize}

  Moreover,
  \begin{itemize}
    \item[(iii)]\begin{equation}\label{eq1-6}
        \int_\C\log^+|z|\,\d\mu_T(z)=
        \frac{1}{2\pi}\int_0^{2\pi}f(\e^{\i\theta})\,\d\theta.
        \end{equation}
  \end{itemize}
\end{thm}

\vspace{.2cm}

The following lemma was proven by F.~Larsen in his unpublished
thesis (cf. \cite[section~2]{FL}. For the convenience of the reader we include a (somewhat
different) proof.

\begin{lemma}\label{Larsen} Let $a,b\in\CM$ and let $\eps>0$. Define
  $g_\eps,g:\C\rightarrow \R$ by
  \[
  g_\eps(\lambda)= \textstyle{\frac12} \tau(\log((a-\lambda
  b)\cc(a-\lambda b) + \eps\unit)),
  \]
  and
  \[
  g(\lambda)=\log\Delta(a-\lambda b).
  \]
  Then $g_\eps$ is subharmonic, and if $g(\lambda)>-\infty$ for some
  $\lambda\in\C$, then $g$ is subharmonic as well.
\end{lemma}

\proof Let $\lambda_1=\re(\lambda)$, $\lambda_2=\im(\lambda)$,
$\lambda\in\C$. At first we show that $(\lambda_1,\lambda_2)\mapsto g_\eps(\lambda_1+\i\lambda_2)$ is a $C^2$--function in
$\R^2$. Fix $\eps>0$, and define $h,k:\C\rightarrow\CM$ by
\begin{eqnarray*}
h(\lambda)&=&(a-\lambda b)\cc (a-\lambda b)+\eps\unit,\\
k(\lambda)&=&(a-\lambda b)(a-\lambda b)\cc + \eps\unit.
\end{eqnarray*}
Then $h$ and $k$ are second order polynomials in
$(\lambda_1,\lambda_2)$ with coefficients in $\CM$, and
$h(\lambda)\geq\eps\unit$, $k(\lambda)\geq \eps\unit$ for all
$\lambda\in\C$. Hence, by \cite[Lemma~4.6]{HT},
\[
g_\eps(\lambda)= \textstyle{\frac12}\tau(\log h(\lambda)), \qquad \lambda\in\C,
\]
has continuous partial derivatives given by
\[
\frac{\partial g_\eps}{\partial \lambda_j}= \textstyle{\frac12}
\tau\Big(h^{-1} \frac{\partial h}{\partial \lambda_j}\Big), \qquad
j=1,2.
\]
Therefore, by \cite[Lemma~3.2]{HT}, $g_\eps$ is a $C^2$--function with
\begin{equation}\label{2.18+}
\frac{\partial^2 g_\eps}{\partial\lambda_i\partial\lambda_j}=
\textstyle{\frac 12}\tau\Big(-h^{-1}\frac{\partial h}{\partial
  \lambda_i}h^{-1}\frac{\partial h}{\partial
  \lambda_j} + h^{-1}\frac{\partial^2
  h}{\partial\lambda_i\partial\lambda_j}\Big),\qquad i=1,2,\; j=1,2.
\end{equation}
Since $g_\eps$ is $C^2$, $g_\eps$ is subharmonic if and only if its
Laplacian
\[
\frac{\partial^2 g_\eps}{\partial\lambda_1^2}+\frac{\partial^2
  g_\eps}{\partial\lambda_2^2}
\]
is positive. Following standard notation, we let
\[
\frac{\partial}{\partial\lambda}=\frac12\Big(\frac{\partial}{\partial\lambda_1}
- \i \frac{\partial}{\partial\lambda_2}\Big) \;\;\; {\rm and} \;\;\;
\frac{\partial}{\partial \overline\lambda}=\frac12\Big(\frac{\partial}{\partial\lambda_1}
+ \i \frac{\partial}{\partial\lambda_2}\Big).
\]
Then
\[
\frac{\partial^2 g_\eps}{\partial\lambda_1^2}+\frac{\partial^2
  g_\eps}{\partial\lambda_2^2} = 4\frac{\partial^2 g_\eps}{\partial\overline{\lambda}\partial\lambda}
\]
By application of \eqref{2.18+}, we find that
\begin{equation}\label{2.18++}
  \frac{\partial^2 g_\eps}{\partial\overline{\lambda}\partial\lambda}=
  \textstyle{\frac 12}\tau\Big(-h^{-1}\frac{\partial h}{\partial
  \overline\lambda}h^{-1}\frac{\partial h}{\partial \lambda} +
  h^{-1}\frac{\partial^2 h}{\partial \overline\lambda
  \partial\lambda}\Big).
\end{equation}
Since
\[
h(\lambda)= a\cc a -\lambda a\cc b -\overline\lambda b\cc a +
|\lambda|^ 2 b\cc b +\eps\unit,
\]
we have
\[
\frac{\partial h}{\partial \lambda} = -a\cc b +\overline\lambda b\cc
b=-(a-\lambda b)\cc b,
\]
\[
\frac{\partial h}{\partial \overline\lambda} = -b\cc a +\lambda b\cc
b=-b\cc(a-\lambda b),
\]
and
\[
\frac{\partial^2 h}{\partial\overline\lambda \partial\lambda}=b\cc b.
\]

Applying the identity $x(x\cc x + \eps\unit)^{-1}= (xx\cc
+\eps\unit)^{-1}x$ to $x= a-\lambda b$, we find that
\begin{eqnarray*}
\frac{\partial^2 h}{\partial\overline\lambda \partial\lambda}-
\frac{\partial h}{\partial\overline\lambda}h^{-1}\frac{\partial
  h}{\partial \lambda} & = &
b\cc b - b\cc x(x\cc x + \eps\unit)^{-1}x\cc b\\
&=& b\cc b - b\cc(xx\cc + \eps\unit)^{-1}xx\cc b\\
&=& b\cc b - b\cc(\unit-\eps(xx\cc + \eps\unit)^{-1})b\\
&=& \eps b\cc(xx\cc + \eps\unit)^{-1}b\\
&=& \eps b\cc k^{-1}b.
\end{eqnarray*}
Then by \eqref{2.18++},
\begin{eqnarray*}
 \frac{\partial^2 g_\eps}{\partial\overline{\lambda}\partial\lambda}&=&
 \textstyle{\frac{\eps}{2}} \tau(h(\lambda)^{-1}b\cc k(\lambda)^{-1}b)\\
 &=& \textstyle{\frac{\eps}{2}} \tau(h(\lambda)^{-\frac12}b\cc
 k(\lambda)^{-1}bh(\lambda)^{-\frac12})\\
 &\geq & 0,
\end{eqnarray*}
showing that $g_\eps$ is subharmonic.

Fix $\lambda\in\C$, and let $x=a-\lambda b$ as above. Then
\[
g_\eps(\lambda)=\textstyle{\frac12}\int_0^{\|x\|}\log(t^2+\eps)\,\d\mu_{|x|}(t),
\]
and
\[
g(\lambda)=\textstyle{\frac12}\int_0^{\|x\|}\log(t^2)\,\d\mu_{|x|}(t).
\]
Hence, $g_\eps$ is a monotonically decreasing function of $\eps>0$,
and
\[
g(\lambda)=\lim_{\eps\rightarrow 0^+}g_\eps(\lambda).
\]
According to \cite{HK}, $g$ is then either subharmonic or identically
$-\infty$. $\endproof$

\vspace{.2cm}

\begin{prop} Let $T\in\CMD$. Then the function $f:\C\rightarrow
  [-\infty,\infty[$ given by
  \[
  f(\lambda)=\log \Delta(T-\lambda\unit)
  \]
  is subharmonic in $\C$.
\end{prop}

\proof Define $T_1$, $T_2\in\CM$ by
\begin{equation}\label{eq1-8}
  T_1 = T(T\cc T+\unit)^{-\frac12}
\end{equation}

and
\begin{equation}\label{eq1-9}
  T_2 = (T\cc T+\unit)^{-\frac12}.
\end{equation}

Then for every $\lambda\in\C$,
\begin{equation*}
  T-\lambda\unit =(T_1-\lambda T_2)T_2^{-1},
\end{equation*}
where $\Delta(T_2)>0$ (cf. \eqref{eq1-11}). Thus, $ T-\lambda\unit\in
\CMD$ with 
\begin{equation*}
  \Delta(T-\lambda\unit)=\Delta(T_1-\lambda T_2)\Delta(T_2)^{-1},
\end{equation*}

i.e.
\begin{equation}\label{eq1-12}
  f(\lambda)=L(T-\lambda\unit)=L(T_1-\lambda T_2)-L(T_2).
\end{equation}

Then by Lemma~\ref{Larsen}, $f$ is either subharmonic or identically
$-\infty$. With
\[
h(\lambda)= L(T_2-\lambda T_1)-L(T_2),
\]
$h(0)=0>-\infty$, and it follows from Lemma~\ref{Larsen} that $h$ is
subharmonic. In particular, $h(\lambda)>-\infty$ for almost every
$\lambda\in\C$ w.r.t. Lebesgue measure. For
$\lambda\in\C\setminus\{0\}$,
\[
f(\lambda)=h\Big(\frac1\lambda\Big) +\log|\lambda|.
\]
Hence, $f$ is not identically $-\infty$. $\endproof$

\vspace{1cm}

Recall from \cite[Section~3.5.4]{HK} that one can associate to every
subharmonic function $u$ the socalled {\it Riesz measure} $\mu_u$, which is a
positive Borel measure on $\R^2$ uniquely determined by
\begin{equation}\label{eq1-30}
  \forall \phi\in C_c^\infty(\R^2):\quad \frac{1}{2\pi}\int_{\R^2} u
  \nabla^2 \phi \,\d m  = \int_{\R^2} \phi\,\d\mu_u.
\end{equation}

One uses the notation $\d\mu_u = \frac{1}{2\pi}\nabla^2u \,\d\lambda$, and
this is what is meant by \eqref{eq1-5}.

\vspace{.2cm}

In order to prove the rest of Theorem~\ref{mainthm}, we need some general lemmas on
subharmonic functions:

\begin{lemma}\label{lemmaB} Let $g:\C\rightarrow [-\infty,\infty[$ be a subharmonic
  function, and for $r>0$ define
  \begin{eqnarray}
    m(g,r)&=&\frac{1}{2\pi} \int_0^{2\pi} g(r\e^{\i\theta})\,\d \theta,\\
    M(g,r)&=& \sup_{|z|=r}g(z).
  \end{eqnarray}
Then
\begin{equation}\label{g(0)}
  g(0)=\lim_{r\rightarrow 0}m(g,r)=\lim_{r\rightarrow 0}M(g,r).
\end{equation}
\end{lemma}

\proof Clearly, $m(g,r)\leq M(g,r)$ for every $r>0$. Moreover, since $g$ is
subharmonic, $g(0)\leq m(g,r)$, $(r>0)$. It follows that
\begin{equation}\label{eq1-15}
  g(0)\leq \left\{\begin{array}{l}
      \limsup_{r\rightarrow 0}m(g,r)\leq  \limsup_{r\rightarrow 0}M(g,r)\\
      \liminf_{r\rightarrow 0}m(g,r)\leq  \liminf_{r\rightarrow 0}M(g,r)
      \end{array}\right .
\end{equation}

Now, every upper semicontinuous function attains a maximum on every compact
set. In particular, there exists for every $r>0$ a complex number $z_r$ of
modulus $r$ such that $g(z_r)=M(g,r)$. $z_r\rightarrow
0$ as $r\rightarrow 0$, and therefore
\begin{equation}\label{eq1-16}
  g(0)\geq \limsup_{r\rightarrow 0}g(z_r) =  \limsup_{r\rightarrow
  0}M(g,r).
\end{equation}

It follows from \eqref{eq1-15} and \eqref{eq1-16} that
\begin{eqnarray*}
  g(0)&\leq& \liminf_{r\rightarrow 0}m(g,r) \\
  & \leq & \left\{\begin{array}{l}
      \limsup_{r\rightarrow 0}m(g,r)\\
      \liminf_{r\rightarrow 0}M(g,r)
      \end{array}\right\}\\
    & \leq &  \limsup_{r\rightarrow 0}M(g,r) \\
    &\leq & g(0),
\end{eqnarray*}
so the four inequalities above are in fact identities, and this proves \eqref{g(0)}. $\endproof$

\vspace{.2cm}

\begin{lemma}\label{lemmaD} $f$ given by \eqref{eq1-4} satisfies
  \begin{equation}
    \lim_{r\rightarrow\infty}(M(f,r)-\log
    r)=\lim_{r\rightarrow\infty}(m(f,r)-\log r) = 0.
  \end{equation}
\end{lemma}

\proof Define $h:\C\rightarrow [-\infty,\infty[$ by
\begin{equation}
 h(\lambda)= L(T_2-\lambda T_1)-L(T_2), \qquad \lambda\in\C.
\end{equation}

Then $h$ is subharmonic with $h(0)=0$, and it follows from
Lemma~\ref{lemmaB} that
\begin{equation}\label{eq1-18}
  0 =  \lim_{r\rightarrow 0}m(h,r) = \lim_{r\rightarrow 0}M(h,r).
\end{equation}

Since
\begin{equation}
  h(\lambda)= \log|\lambda| + f(\textstyle\frac{1}{\lambda}), \qquad
  \lambda\neq0,
\end{equation}

we get that when $r>0$,
\begin{eqnarray*}
  M(f,r)&=&M(h, \textstyle\frac 1r) + \log r,\\
  m(f,r)&=& m(h, \textstyle\frac 1r) + \log r,
\end{eqnarray*}

and combining this with \eqref{eq1-18} we obtain the desired result. $\endproof$

\vspace{.2cm}

\begin{lemma}\label{lemmaC} Let $R>r>0$, and let $g$ be subharmonic in $\C$. Then with $\d\mu
    = \frac{1}{2\pi}\nabla^2 g \,\d\lambda$ and
    \[
    \psi(z)=\left\{\begin{array}{lll}
        \log\big(\frac Rr\big) &,& |z|\leq r\\
        \log\big(\frac{R}{|z|}\big) &, & r<|z|<R\\
        0 &,& |z|\geq R
        \end{array}\right.
    \]
  one has that
  \begin{equation}
    m(g,R)-m(g,r)=\int_\C \psi(z)\,\d\mu(z).
  \end{equation}
\end{lemma}

\proof Cf. \cite[(3.5.7)]{HK}. $\endproof$

\vspace{.2cm}

{\it Proof of Theorem~\ref{mainthm}.} When $R>1>0$ define
$\psi_R:\C\rightarrow \R$ by
\[
    \psi_R(z)=\left\{\begin{array}{lll}
        \log R &,& |z|\leq 1\\
        \log\big(\frac{R}{|z|}\big) &, & 1<|z|<R\\
        0 &,& |z|\geq R
        \end{array}\right.
    \]

Then, according to Lemma~\ref{lemmaC},
\begin{equation}\label{eq1-19}
  \int_\C \psi_R(z)\,\d\mu_T(z) = m(f,R)-m(f,1).
\end{equation}

Now, $\frac{1}{\log R}\psi_R \nearrow 1$ as $R\rightarrow \infty$, so by
the Monotone Convergence Theorem, \eqref{eq1-19} and Lemma~\ref{lemmaD},
\[
\mu_T(\C)=\lim_{R\rightarrow\infty}\,\frac{m(f,R)-m(f,1)}{\log R} =
1,
\]

that is, $\mu_T$ is a probability measure.

When $R>1$, let
\begin{equation}
  \omega_R(z)=\log R -\psi_R(z), \qquad z\in\C.
\end{equation}

Then $\omega_R(z)\nearrow \log^+|z|$ as $\R\rightarrow \infty$, and hence
by one more application of Lemma~\ref{lemmaD},
\begin{eqnarray*}
  \int_\C\log^+|z|\,\d\mu_T(z)&=& \lim_{R\rightarrow\infty}
  \int_\C\omega_R\,\d\mu_T\\
  &=&  \lim_{R\rightarrow\infty}(\log R - m(f,R)+m(f,1))\\
  &=& m(f,1),
\end{eqnarray*}

proving \eqref{eq1-6}. Note that since $f$ is subharmonic,  \eqref{eq1-6}
imlies that $ \int_\C\log^+|z|\,\d\mu_T(z)<\infty$.

To see that \eqref{eq1-7} holds, it suffices to consider the case
$\lambda=0$. Indeed, for fixed $\lambda\in\C$ one easily sees that
$\mu_{T-\lambda\unit}$ is the push-forward measure of $\mu_T$ under the map
$z\mapsto z-\lambda$ (cf. Lemma~\ref{lemmaJ}), and therefore
\begin{equation}
\int_\C\log|z-\lambda|\,\d\mu_T(z)=
\int_\C\log|z|\,\d\mu_{T-\lambda\unit}(z).
\end{equation}

In the case $\lambda =0$ one has to compute the integrals
$\int_\C\log^\pm |z|\,\d\mu_T(z)$. We have just seen that
\begin{equation}\label{eq1-20}
\int_\C\log^+|z|\,\d\mu_T(z)=m(f,1),
\end{equation}

and with
\[
\chi_r(z)=\left\{\begin{array}{lll}
      \log\frac 1r &,& |z|\leq r\\
      \log\frac{1}{|z|} &,& r<|z|\leq 1\\
      0 &,& |z|\geq 1
      \end{array}
      \right.
\]

$\chi_r(z)\nearrow \log^-|z|$ as $r\searrow 0$. Hence by Lemma~\ref{lemmaB}
and Lemma~\ref{lemmaC},
\begin{eqnarray*}
  \int_\C\log^- |z|\,\d\mu_T(z) & = & \lim_{r\rightarrow
  0}\int_\C\chi_r\,\d\mu_T\\
& = & \lim_{r\rightarrow 0} (m(f,1)-m(f,r))\\
& = & m(f,1)-f(0).
\end{eqnarray*}

Combining this with \eqref{eq1-20} we get that
\[
\int_\C\log|z|\,\d\mu_T(z) = f(0)=L(T),
\]

as desired. 

In order to prove that $\mu_T$ is uniquely determined by (i) and (ii) of
Theorem~\ref{mainthm}, suppose $\nu\in {\rm Prob}(\C)$ satisfies
\begin{equation}\label{firstcond}
  \int_\C \log^+|z|\,\d\nu(z)<\infty,
\end{equation}
and
\begin{equation}
  \forall\; \lambda\in\C:\qquad \int_\C\log|z-\lambda|\,\d\nu(z) =
  L(T-\lambda\unit).
\end{equation}

Note that \eqref{firstcond} implies that $\int_\C\log|z-\lambda|\,\d\nu(z)$
is well-defined, since
\[
\log|z-\lambda|\leq \log(|z|+|\lambda|),
\]
and
\[
|z|+|\lambda|\leq (|\lambda + 1|)\cdot \max\{1, |z|\}.
\]
Hence
\begin{equation}\label{logineq}
  \log|z-\lambda|\leq \log(|\lambda|+1)+ \log^+|z|.
\end{equation}

Since $\mu$ and $\nu$ are both probability measures, it follows from a
$C^\infty$-version of Urysohn's Lemma (cf. \cite[(8.18)]{Fo}) that if
\[
\int_\C \phi\,\d\mu_T = \int_\C \phi\,\d\nu
\]
for every function $\phi\in C_c^\infty(\R^2)$, then $\mu_T = \nu$. Then
consider an arbitrary function $\phi\in C_c^\infty(\R^2)$. Since the
Laplacian of $w\mapsto \frac{1}{2\pi}\log|w-z|$ (in the distribution sense)
is the Dirac measure $\delta_z$ at $z$, one has that
\begin{eqnarray}\label{noFubini}
  \int_\C\phi(z)\,\d\nu(z)&=& \int_\C \Big(\int_\C
  \phi(\lambda)\,\delta_z(\lambda)\Big)\,\d\nu(z)\nonumber\\
  &=& \frac{1}{2\pi} \int_\C\Big(\int_\C
  (\nabla^2\phi)(\lambda)\log|z-\lambda|\,\d\lambda\Big) \,\d\nu(z).
\end{eqnarray}

At this place we would like to reverse the order of integration, but
it is not entirely clear that this is a legal operation. Therefore we put $M=\|\nabla^2\phi\|_\infty$, and take $\chi\in C_c^\infty(\R^2)$ such
that $0\leq \chi\leq 1$ and $\chi_{|_{\supp(\nabla^2\phi)}} = 1$. With
\[
\psi_1 = \frac12(M+\nabla^2\phi)\chi
\]
and
\[
\psi_2 = \frac 12(M-\nabla^2\phi)\chi
\]
one has that $\psi_1, \psi_2\in C_c^\infty(\R^2)^+$, and
$\nabla^2\phi = \psi_1 -\psi_2$. 

Also not that, according to \eqref{logineq},
\[
h(\lambda,z):= \log(|\lambda|+1)+\log^+|z| - \log|z-\lambda|\geq 0.
\]

Therefore by Tonelli's Theorem
\begin{equation}\label{Tonelli}
\int_\C \psi_i(\lambda)\int_\C h(\lambda,z)\,\d\nu(z) \,\d\lambda =
\int_\C\int_\C \psi_i(\lambda) h(\lambda, z)\,\d\lambda \,\d\nu(z), \qquad i=1,2.
\end{equation}

The map $\lambda \mapsto L(T-\lambda\unit)$ is subharmonic and therefore
locally integrable.
Since
\[
\int_\C h(\lambda,z)\,\d\nu(z)= \log(|\lambda|+1)+
\int_\C\log^+|z|\,\d\nu(z) - L(T-\lambda\unit),
\]
where $\lambda \mapsto L(T-\lambda\unit)$ is subharmonic and therefore
locally integrable,
\[
\int_\C \psi_i(\lambda)\int_\C h(\lambda,z)\,\d\nu(z) \,\d\lambda <\infty,
\qquad i=1,2.
\]
It now follows from \eqref{Tonelli} that
\[
\int_\C(\nabla^2\phi)(\lambda)\int_\C h(\lambda,z)\,\d\nu(z) \,\d\lambda =
\int_\C\int_\C (\nabla^2\phi)(\lambda) h(\lambda, z)\,\d\lambda \,\d\nu(z),
\]
and since
\[
\int_\C|(\nabla^2\phi)(\lambda)|\int_\C \log(|\lambda|+1)\,\d\nu(z)
\,\d\lambda <\infty,
\]
and
\[
\int_\C|(\nabla^2\phi)(\lambda)|\int_\C \log^+|z|\,\d\nu(z)
\,\d\lambda <\infty,
\]
we deduce that
\begin{eqnarray*}
\int_\C\phi(z)\,\d\nu(z) & = & \frac{1}{2\pi}\int_\C\Bigg(\int_\C (\nabla^2\phi)(\lambda) \log|\lambda-z|
\,\d\lambda\Bigg) \,\d\nu(z) \\
&=&
\frac{1}{2\pi}\int_\C(\nabla^2\phi)(\lambda)\int_\C
\log|\lambda-z|\,\d\nu(z) \,\d\lambda \\
& = & \frac{1}{2\pi}\int_\C (\nabla^2\phi)(\lambda) L(T-\lambda\unit)
\,\d\lambda \\
& = & \int_\C \phi(z)\,\d\mu_T(z),
\end{eqnarray*}
and this is the desired identity.  $\endproof$

\vspace{.2cm}

It follows from Theorem~\ref{mainthm} that one can associate to every
operator $T\in\CMD$ a probability measure $\mu_T$ on $(\C,\B_2)$, such that
in the case where $T\in\CM$,  $\mu_T$ agrees with the Brown measure of
$T$. Therefore we make the following definition:

\begin{definition} For $T\in\CMD$ we shall say that the probability measure
  $\mu_T$ from Theorem~\ref{mainthm} is the {\it Brown measure} of
  $T$. 
\end{definition}

\vspace{.2cm}

In the remaining part of this section we will see that many of the properties of
the Brown measure for bounded operators carry over to this more general
setting.

\begin{prop}\label{lemmaJ}Let $T\in\CMD$. Then for every $r>0$ and every $\lambda\in\C$, the Brown measure of $rT+\lambda\unit$,
  $\mu_{rT+\lambda\unit}$, is the push-forward measure of $\mu_T$ via the
  map $z\mapsto rz + \lambda$.
\end{prop}

\proof Making use of Urysohn's Lemma for $C^\infty$-functions on $\R^2$
(cf. \cite[(8.18)]{Fo}) and the fact that both of the measures considered
here are probability measures, one easily sees that if  
\[
\int_\C \phi(z)\,\d\mu_{rT+\lambda\unit}(z)=\int_\C \phi(rz+\lambda)\,\d\mu_T(z)
\]
for every $\phi\in C_c^\infty(\R^2)$, then the two measures in speak agree on
compact sets and hence on all of $\B_2$.

Let $\phi\in C_c^\infty(\R^2)$. Then by definition,
\begin{eqnarray*}
\int_\C \phi(rz+\lambda)\,\d\mu_T(z) & = & \frac{1}{2\pi} \int_\C
\Big(\frac{\partial^2}{\partial z_1^2}+ \frac{\partial^2}{\partial
  z_2^2}\Big)\phi(rz+\lambda)\,f(z)\,\d z\\
&=& \frac{1}{2\pi} \int_\C
r^2\Big(\frac{\partial^2}{\partial w_1^2}+ \frac{\partial^2}{\partial
  w_2^2}\Big)\phi(w)\,f\Big(\frac1r(w-\lambda)\Big)\frac{1}{r^2}\,\d w\\
& = & \frac{1}{2\pi} \int_\C
\nabla^2\phi(w)f\Big(\frac1r(w-\lambda)\Big)\,\d w\\
& = &  \frac{1}{2\pi} \int_\C
\nabla^2\phi(w)\,[L(rT+\lambda\unit -w\unit)-\log r]\,\d w\\
& = & \int_\C \phi(w)\,\d\mu_{rT+\lambda\unit}(w) - \log r \cdot \int_\C
\nabla^2\phi(w)\,\d w\\
& = & \int_\C \phi(w)\,\d\mu_{rT+\lambda\unit}(w),
\end{eqnarray*}

where the last identity follows from Green's Theorem. $\endproof$

\vspace{.2cm}

\begin{prop}\label{power}For every $T\in\CMD$ and every $m\in\N$,
  $\mu_{T^m}$ is the
  push-forward measure of $\mu_T$ via the map $z\mapsto z^m$.
\end{prop}

\proof Let $\nu\in {\rm Prob}(\C)$ denote the push-forward measure of
$\mu_T$ under the map $z\mapsto z^m$. According to Theorem~\ref{mainthm} it
suffices to prove that
\[
\int_\C\log^+|z|\,\d\nu(z)<\infty,
\]
and
\[
\forall\,\lambda\in\C: \qquad \int_\C \log|\lambda-z|\,\d\nu(z) =
L(T^m-\lambda\unit).
\]

Here
\begin{eqnarray*}
\int_\C\log^+|z|\,\d\nu(z)&=&\int_\C\log^+|z^m|\,\d\mu_T(z) \\
&=& m\int_\C\log^+|z|\,\d\mu_T(z)\\
&<&\infty,
\end{eqnarray*}
and if we let $\theta_1,\ldots, \theta_m$ denote the $m$ complex roots of
$Q(z)=z^m-1$, then for every $\lambda\in\C$,
\[
|\lambda-z^m|=\prod_{k=1}^m|\theta_k\lambda^{\frac1m}-z|.
\]
Hence
\begin{eqnarray*}
 \int_\C \log|\lambda-z|\,\d\nu(z)& = & \int_\C
 \log|\lambda-z^m|\,\d\mu_T(z)\\
 & = & \int_\C\sum_{k=1}^m\log|\theta_k\lambda^{\frac1m}-z|\,\d\mu_T(z)\\
 & = & \sum_{k=1}^m L(T-\theta_k\lambda^{\frac1m}\unit)\\
 & = & L\Big(\prod_{k=1}^m (T-\theta_k\lambda^{\frac1m}\unit)\Big)\\
 & = & L(T^m-\lambda\unit),
\end{eqnarray*}
as desired. $\endproof$

\vspace{.2cm}

\begin{prop}\label{inverse} If $T\in\CMD$ with
  \begin{equation}\label{app3}
    \int_0^1 \log t\,\d\mu_{|T|}(t) >-\infty,
  \end{equation}
then $\mu_T(\{0\})=\mu_{|T|}(\{0\})=0$, and $T$ has an inverse
$T^{-1}\in\CMD$. Moreover, $\mu_{T^{-1}}$ is the push-forward 
measure of $\mu_T$ via the map $z\mapsto z^{-1}$.
\end{prop}

\proof According to Theorem~\ref{mainthm},
\begin{equation}
  \int_\C \log|z|\,\d\mu_T(z)= L(T)=\int_0^\infty \log t\,\d\mu_{|T|}(t).
\end{equation}
Hence, if \eqref{app3} holds, then
\begin{equation}\label{app4}
-\infty <\int_\C \log|z|\,\d\mu_T(z)<\infty,
\end{equation}
and therefore $\mu_T(\{0\})=\mu_{|T|}(\{0\})=0$. Moreover, $|T|$ has an
inverse $|T|^{-1}\in \tilde\CM$ with
\begin{eqnarray*}
  \int_0^\infty \log^+(t)\,\d\mu_{|T|^{-1}}(t) &=& \int_0^\infty
  \log^+\Big(\frac1t\Big)\,\d\mu_{|T|}(t)\\
  &=& -\int_0^1\log t\,\d\mu_{|T|}(t)\\
  & <&\infty,
\end{eqnarray*}
so $|T|^{-1}\in\CMD$. Take $U\in\CU(\CM)$ such that $T=U|T|$. Then
$T^{-1}=|T|^{-1}U\cc\in\CMD$.

Now, let $\nu$ denote the push-forward measure of $\mu_{T}$ under the map 
$z\mapsto z^{-1}$. According to Theorem~\ref{mainthm}, if
\begin{equation}
  \int_\C \log^+|z|\,\d\nu(z) <\infty,
\end{equation}
and
\begin{equation}\label{app5}
 \forall\,\lambda\in\C: \qquad \int_\C\log|z-\lambda|\,\d\nu(z) = L(T^{-1}-\lambda\unit), 
\end{equation}
then $\nu = \mu_{T^{-1}}$. Applying \eqref{app4} we find that
\begin{eqnarray*}
   \int_\C \log^+|z|\,\d\nu(z) &=& \int_\C
   \log^+\Big|\frac1z\Big|\,\d\mu_T(z)\\
   &=& -\int_{(|z|\leq 1)} \log|z|\,\d\mu_T(z)\\
   &<&\infty.
\end{eqnarray*}
In order to prove that \eqref{app5} holds, let $\lambda\in\C$. If
$\lambda\neq 0$, then, using the multiplicativity of $\Delta$ on $\CMD$, we
find that
\begin{eqnarray*}
  \int_\C\log|z-\lambda|\,\d\nu(z) & = &
  \int_\C\log\Big|\frac1z-\lambda\Big|\,\d\mu_T(z)\\
  & = &
  \int_\C\log\Big|\frac1z\Big(\frac1\lambda-z\Big)\lambda\Big|\,\d\mu_T(z)\\
  & = & \int_\C \Big(\log|\lambda| + \log\Big|\frac1\lambda-z\Big| -
  \log|z|\Big)\,\d\mu_T(z)\\
  & = & L(\lambda\unit)+L\Big(T-\frac1\lambda\unit\Big)-L(T)\\
  &=& L\Big(\lambda\unit\Big(T-\frac1\lambda\unit\Big)T^{-1}\Big)\\
  & = & L(T^{-1}-\lambda\unit).
\end{eqnarray*}
In the case $\lambda =0$ we have:
\begin{eqnarray*}
   \int_\C\log|z|\,\d\nu(z) & = & -\int_\C\log|z|\,\d\mu_T(z)\\
   & = & -L(T)\\
   &= & L(T^{-1}).
\end{eqnarray*}
Hence \eqref{app5} holds, and $\nu= \mu_{T^{-1}}$. $\endproof$

\vspace{.2cm}

\begin{prop}
Let $T\in\CMD$. Then $\supp(\mu_T)\subseteq\sigma(T)$.
\end{prop}

\proof Let $\lambda\in\C\setminus\sigma(T)$. Then $T-\lambda\unit$ is
invertible with bounded inverse. Moreover, according to
Proposition~\ref{inverse}, $\mu_{(T-\lambda\unit)^{-1}}$ is the push--forward
  measure of $\mu_{T-\lambda\unit}$ via the map $z\mapsto z^{-1}$,
  $z\in\C\setminus\{0\}$. Since $(T-\lambda\unit)^{-1}$ is bounded, we
  have from \cite{Bro} that
  \[
  \supp(\mu_{(T-\lambda\unit)^{-1}})\subseteq
  \sigma((T-\lambda\unit)^{-1})\subseteq \overline{B(0,r)},
  \]
  where $r=\|(T-\lambda\unit)^{-1}\|$. Hence,
  \[
  \supp(\mu_{T-\lambda\unit})\subseteq \{z\in\C\,|\,|z|\geq
  \textstyle{\frac1r}\}.
  \]
  In particular, $0\notin \supp(\mu_{T-\lambda\unit})$, which by
  Proposition~\ref{lemmaJ} is equivalent to
  $\lambda\notin\supp(\mu_T)$. Hence, $\supp(\mu_T)\subseteq
  \sigma(T)$. $\endproof$

\vspace{.2cm}

\begin{lemma} \label{p}
  For every $p\in(0,\infty)$ and every $t\in [0,\infty[$,
  \begin{equation}\label{app6}
    t^p= p^2 \int_0^\infty \log^+(at)a^{-p-1}\,\d a.
  \end{equation}
\end{lemma}

\proof For $t=0$ this is trivial. For $t>0$ we find that
\begin{eqnarray*}
  \int_0^\infty \log^+(at)a^{-p-1}\,\d a & = & \int_{\frac1t}^\infty
  \log(at)a^{-p-1}\,\d a\\
   & = & \Bigg[-\frac1p \log(at)\,a^{-p} \Bigg]_{\frac1t}^\infty -
  \int_{\frac1t}^\infty -\frac{1}{p\,a}\, a^{-p}\,\d a\\
  & = & 0 - \Bigg[-\frac{1}{p^2}a^{-p} \Bigg]_{\frac1t}^\infty\\
  & = & \frac{1}{p^2} t^p. \quad \endproof
\end{eqnarray*}

\vspace{.2cm}

We will now prove Weil's inequality for operators $T$ in
$L^p(\CM)$ (cf. \cite[corollary~3.8]{Bro} for the case $T\in\CM$):

\begin{thm}\label{propE} Let  $p\in (0,\infty)$ and let $T\in
  L^p(\CM)$. Then
  \begin{equation}
     \int_\C|z|^p\,\d\mu_T(z)\leq \|T\|_p^p.
  \end{equation}
\end{thm}

In the proof of this theorem we shall need the following lemma, the
proof of which we postpone for a while:

\begin{lemma}\label{app7}
 Let $T\in\CMD$. Then
 \begin{equation}
   \int_\C \log^+|z|\,\d\mu_T(z)\leq \tau(\log^+|T|).
 \end{equation}
\end{lemma}

{\it Proof of Proposition~\ref{propE}.} Let $a\geq 0$. Then, according to Lemma~\ref{lemmaJ} and
Lemma~\ref{app7},
\begin{eqnarray*}
  \int_\C\log^+(a|z|)\,\d\mu_T(z)&=& \int_\C\log^+|z|\,\d\mu_{aT}(z)\\
  & \leq & \int_0^\infty \log^+ t\,\d\mu_{|aT|}(t)\\
  & = & \int_0^\infty \log^+(at) \d\mu_{|T|}(t).
\end{eqnarray*}
Hence by Lemma~\ref{p} and Tonelli's Theorem,
\begin{eqnarray*}
  \int_\C |z|^p\d\mu_T(z) & = & p^2\int_0^\infty
  \Bigg(\int_\C\log^+(a|z|)\,\d\mu_T(z)\Bigg)a^{-p-1}\,\d a\\
  & \leq & p^2 \int_0^\infty
  \Bigg(\int_0^\infty \log^+(at)\,\d\mu_{|T|}(t)\Bigg)a^{-p-1}\,\d a\\
  & = & \int_0^\infty t^p \,\d\mu_{|T|}(t)\\
  & = & \tau(|T|^p). \qquad \endproof
\end{eqnarray*}

\vspace{.2cm}

In order to prove Lemma~\ref{app7} we shall need some additional results:

\begin{lemma}\label{lemmaG} Suppose $A, B, C\in\CMD$ with $A$ and $B$
  invertible in $\CMD$ and
  \[
  \begin{pmatrix}
    A & C\cc\\
    C & B
  \end{pmatrix} \geq 0.
  \]
  Then
  \begin{equation}
    \Delta(C)\leq \Delta(A)^{\frac 12}\Delta(B)^{\frac 12}.
  \end{equation}
\end{lemma}

\proof Note that $A, B\geq 0$ and that
\[
 \begin{pmatrix}
    \unit & A^{-\frac 12}C\cc B^{-\frac 12}\\
    B^{-\frac 12}CA^{-\frac 12} & \unit
 \end{pmatrix} =
 \begin{pmatrix}
    A^{-\frac 12} & 0 \\
    0 &  B^{-\frac 12}
 \end{pmatrix}
 \begin{pmatrix}
    A & C\cc\\
    C & B
 \end{pmatrix}
 \begin{pmatrix}
    A^{-\frac 12} & 0 \\
    0 &  B^{-\frac 12}
 \end{pmatrix}\geq 0,
\]

which is equivalent to saying that $\|B^{-\frac 12}CA^{-\frac 12}\|\leq 1$,
and this clearly implies that
\[
\Delta(B^{-\frac 12}CA^{-\frac 12})\leq 1. \qquad \endproof
\]

\begin{lemma}\label{lemmaH} For every $S\in\CMD$,
  \begin{equation}
    \Delta(\unit + S)\leq \Delta(\unit + |S|).
  \end{equation}
\end{lemma}

\proof Take a unitary $U\in\CM$ such that $S=U|S|$. Then
\[
\begin{pmatrix}
  |S| &  |S|\\
  |S| &  |S|
\end{pmatrix}\geq 0,
\]
and
\[
\begin{pmatrix}
  \unit &  U\cc\\
  U &  \unit
\end{pmatrix}\geq 0,
\]
whence
\[
\begin{pmatrix}
  |S|+ \unit & |S|+ U\cc\\
  |S| + U &  |S|+\unit
\end{pmatrix}\geq 0.
\]

Now Lemma~\ref{lemmaG} implies that
\begin{eqnarray*}
\Delta(S+\unit)&=&\Delta(U\cc(S+\unit))\\
&=&\Delta(U\cc(U|S|+\unit))\\
&=&\Delta(|S|+U\cc)\\
&\leq& \Delta(|S|+\unit)^\frac12 \Delta(|S|+\unit)^\frac12 \\
&= &\Delta(|S|+\unit),
\end{eqnarray*}

as desired. $\endproof$

\vspace{.2cm}

\begin{lemma}\label{lemmaI} Every $S\in\CMD$ satisfies
  \begin{equation}
    \Delta(\unit +|S^2|)\leq \Delta(\unit + |S|^2),
  \end{equation}
  implying that for arbitrary $n\in\N$,
  \begin{equation}
    \Delta(\unit +|S^{2^n}|)\leq \Delta(\unit + |S|^{2^n}).
  \end{equation}
\end{lemma}

\proof Take a unitary $U\in\CM$ such that $S^2= U|S^2|$. Since
\[
\begin{pmatrix}
  SS\cc &  S^2\\
  (S\cc)^2 &  S\cc S
\end{pmatrix} =
\begin{pmatrix}
S\\
S\cc
\end{pmatrix}
\begin{pmatrix}
S\cc & S
\end{pmatrix}
\geq 0,
\]

we find as in the foregoing proof that
\[
\begin{pmatrix}
  \unit + SS\cc &  U\cc+ S^2\\
  U + (S\cc)^2 &  \unit + S\cc S
\end{pmatrix}\geq 0.
\]

Again this implies that
\[
\Delta(\unit + |S^2|)= \Delta(S^2+U\cc)\leq \Delta(\unit +S\cc S)^\frac12
\Delta(\unit +SS\cc)^\frac12 = \Delta(\unit +S\cc S), 
\]

where the last identity follows from the fact that $S\cc S$ and $SS\cc$ have the
same distribution w.r.t. $\tau$. $\endproof$

\vspace{.2cm}



\vspace{.2cm}

{\it Proof of Lemma~\ref{app7}.}  According to \eqref{eq1-20} we
have:
\begin{equation}\label{eq1-21}
  \int_\C\log^+|z|\,\d\mu_T(z) =  \frac{1}{2\pi}
  \int_0^{2\pi}f(\e^{\i\theta})\,\d\theta ,
\end{equation}

where
\begin{equation}
  f(\lambda)= \tau(\log|T-\lambda\unit|) = \log\Delta(T-\lambda\unit), \qquad \lambda\in\C.
\end{equation}

For every positive integer $n$ define $f_{n}$ by
\begin{equation}
  f_{n}(z)= \sum_{k=0}^{2^n-1} f\big(\e^{\frac{2\pi k}{2^n}\i}z\big),\qquad
  z\in\C.
\end{equation}

Then clearly,
\begin{equation}\label{eq1-22}
\frac{1}{2\pi} \int_0^{2\pi}f(\e^{\i\theta})\,\d\theta =
\frac{1}{2\pi2^n} \int_0^{2\pi}f_{n}(\e^{\i\theta})\,\d\theta .
\end{equation}

Applying Lemma~\ref{lemmaH} and Lemma~\ref{lemmaI} we obtain an estimate of $ f_{n}(\e^{\i\theta})$:
\begin{eqnarray*}
  f_{n}(\e^{\i\theta}) & = & \sum_{k=0}^{2^n-1}\log
  \Delta\big(\e^{-\i\theta}\e^{-\frac{2\pi k}{2^n}\i}T-\unit\big)\\
  &= & \log \Delta\Big(\prod_{k=0}^{2^n-1}\big(\e^{-\i\theta}\e^{-\frac{2\pi
  k}{2^n}\i}T-\unit\big)\Big)\\
  & = & \log \Delta\big(\unit - \e^{-\i2^n \theta}T^{2^n}\big)\\
  & \leq & \log \Delta\big(\unit + |T^{2^n}|\big)\\
  & \leq & \log \Delta\big(\unit + |T|^{2^n}\big)\\
  & = & \tau\big(\log\big(\unit + |T|^{2^n}\big)\big).
\end{eqnarray*}

Combining \eqref{eq1-21} and \eqref{eq1-22} with the above estimate we
see that
\begin{eqnarray*}
 \int_\C\log^+|z|\,\d\mu_T(z) & \leq & \textstyle{\frac{1}{2^n}}\tau\big(\log\big(\unit +
 |T|^{2^n}\big)\big)\\
 & = &  \frac{1}{2^n} \int_{[0,\infty[}
 \log(1+t^{2^n})\,\d\mu_{|T|}(t)\\
 & \leq &  \frac{1}{2^n} \int_{[0,1[}\log 2 \,\d\mu_{|T|}(t) +
 \frac{1}{2^n}  \int_{[1,\infty[} (\log 2+2^n\log t)
 \,\d\mu_{|T|}(t)\\
 & \leq & \frac{2\log 2}{2^n} + \int_{[0,\infty[} \log^+t\,\d\mu_{|T|}(t).
\end{eqnarray*}

Finally, let $n\rightarrow \infty$, and conclude that
\begin{equation*}
 \int_\C\log^+|z|\,\d\mu_T(z) \leq
\int_{[0,\infty[} \log^+t\,\d\mu_{|T|}(t). \qquad \endproof
\end{equation*}

\vspace{.2cm}

\vspace{.2cm}






\vspace{.2cm}

\begin{prop} Let $T\in\CMD$, and suppose $P\in\CM$
  is a non-trivial $T$-invariant projection, i.e. $PTP=TP$. Then
  \begin{equation}\label{eq1-27}
    \Delta(T)=\Delta_{P\CM P}(PTP)^{\tau(P)}\Delta_{P^\bot\CM
    P^\bot}(P^\bot TP^\bot)^{1-\tau(P)},
  \end{equation}
where $\Delta_{P\CM P}$ and $\Delta_{P^\bot\CM P^\bot}$ refer to the
Fuglede-Kadison determinant computed relative to the normalized traces
$\frac{1}{\tau(P)}\tau|_{P\CM P}$ and
$\frac{1}{\tau(P^\bot)}\tau|_{P^\bot \CM P^\bot}$ on $P\CM
P$ and $P^\bot \CM P^\bot$, respectively.
\end{prop}

\proof Put $T_{11}=PTP$, $T_{12}= PTP^\bot$ and $T_{22}=P^\bot TP^\bot$. Then, w.r.t. to the decomposition $\CH = P(\CH)\oplus P(\CH)^\bot$, we
may write
\[
T=\begin{pmatrix}
  T_{11} & T_{12}\\
  0 & T_{22}
  \end{pmatrix}
  =
  \begin{pmatrix}
  \unit & 0\\
  0 & T_{22}
  \end{pmatrix}
  \begin{pmatrix}
  \unit & T_{12}\\
  0 & \unit
  \end{pmatrix}
  \begin{pmatrix}
  T_{11} & 0\\
  0 & \unit
  \end{pmatrix},
\]
where
\[
\Delta\Bigg( \begin{pmatrix}
  \unit & 0\\
  0 & T_{22}
  \end{pmatrix}\Bigg) = \Delta_{P^\bot\CM
    P^\bot}(P^\bot TP^\bot)^{1-\tau(P)},
\]
and
\[
\Delta\Bigg( \begin{pmatrix}
  T_{11} & 0\\
  0 & \unit
  \end{pmatrix}\Bigg) = \Delta_{P\CM
    P}(P TP)^{\tau(P)}.
\]
Thus, \eqref{eq1-27} holds if
\begin{equation}
  \Delta\Bigg( \begin{pmatrix}
  \unit & T_{12}\\
  0 & \unit
  \end{pmatrix}\Bigg)=1.
\end{equation}
To that \eqref{eq1-27} holds, note that
\[
 \begin{pmatrix}
  \unit & T_{12}\\
  0 & \unit
  \end{pmatrix}^{-1} =  \begin{pmatrix}
  \unit & -T_{12}\\
  0 & \unit
  \end{pmatrix},
\]
and hence
\begin{equation}\label{determinant}
\Delta\Bigg( \begin{pmatrix}
  \unit & T_{12}\\
  0 & \unit
  \end{pmatrix}\Bigg)\Delta\Bigg( \begin{pmatrix}
  \unit & -T_{12}\\
  0 & \unit
  \end{pmatrix}\Bigg)=1.
\end{equation}
Also,
\[
\begin{pmatrix}
  \unit & -T_{12}\\
  0 & \unit
  \end{pmatrix} = \begin{pmatrix}
  \unit & 0\\
  0 & -\unit
  \end{pmatrix}\begin{pmatrix}
  \unit & T_{12}\\
  0 & \unit
  \end{pmatrix}\begin{pmatrix}
  \unit & 0\\
  0 & -\unit
  \end{pmatrix},
\]
so that
\[
\Delta\Bigg(\begin{pmatrix}
  \unit & -T_{12}\\
  0 & \unit
  \end{pmatrix}\Bigg)=\Delta\Bigg(\begin{pmatrix}
  \unit & T_{12}\\
  0 & \unit
  \end{pmatrix}\Bigg),
\]
and then by \eqref{determinant},
\[
\Delta\Bigg(\begin{pmatrix}
  \unit & T_{12}\\
  0 & \unit
  \end{pmatrix}\Bigg)=1,
\]
as desired. $\endproof$

\vspace{.2cm}

\begin{lemma}\label{Lp-cont} Let $p\in (0,\infty)$, and let $\eps>0$. Then the map $L_\eps:
  \Lp\rightarrow \R$ given by
  \begin{equation}
    L_\eps(T)= \textstyle{\frac12}\, \tau(\log(T\cc T + \eps\unit)), \qquad T\in \Lp,
  \end{equation}
is continuous w.r.t. $\|\cdot\|_p$.
\end{lemma}

\proof Suppose $T,$ $T_n\in \Lp$ with
\[
\lim_{n\rightarrow \infty}\|T-T_n\|_p =0.
\]

Then $\lim_{n\rightarrow \infty}\|T\cc T-T_n\cc T_n\|_{\frac p2} =0$,
implying that $T_n\cc T_n \rightarrow T\cc T$ in the measure
topology. Consequently,
\begin{equation}\label{eq1-23}
\mu_{T\cc T} = w\cc-\lim_{n\rightarrow\infty}\mu_{T_n\cc T_n}.
\end{equation}

Define a sequence $(\nu_n)_{n=1}^\infty$ of (finite) measures on $(\R, \B)$
by
\begin{equation}
  \d \nu_n(t)=(1+t^{\frac p2})\d\mu_{T_n\cc T_n}(t),
\end{equation}

and note that  since $\lim_{n\rightarrow\infty}\|T_n\|_p = \|T\|_p$,
\begin{equation}\label{eq1-24}
  \sup_{n\in\N}\nu_n(\R)<\infty.
\end{equation}

Similarly define a finite measure $\nu$ on  $(\R, \B)$ by
\begin{equation}
  \d\nu(t)= (1+t^{\frac p2})\d\mu_{T\cc T}(t).
\end{equation}

Because of \eqref{eq1-23} we have that for every $\phi\in C_c(\R)$,
\begin{equation}\label{eq1-25} 
  \int_0^\infty \phi(t)\d\nu(t) = \lim_{n\rightarrow\infty} \int_0^\infty
  \phi(t)\d\nu_n(t).
\end{equation}

When $\phi\in C_0(\R)$, $\phi$ may be approximated (uniformly) by functions
from $C_c(\R)$. Thus, taking \eqref{eq1-24} and \eqref{eq1-25} into account,
one easily sees that
\begin{equation}\label{eq1-26}
  \int_0^\infty \phi(t)\d\nu(t) = \lim_{n\rightarrow\infty} \int_0^\infty
  \phi(t)\d\nu_n(t).
\end{equation}

In particular, with
\begin{equation}
  \phi(t)=\frac{\log(t+\eps)}{1+t^{\frac p2}}, \qquad (t\geq 0),
\end{equation}

\eqref{eq1-26} implies that
\[
L_\eps(T) = \int_0^\infty \phi(t)\d\nu(t) = \lim_{n\rightarrow\infty}
\int_0^\infty \phi(t)\d\nu_n(t)=\lim_{n\rightarrow\infty}L_\eps(T_n). \quad
\endproof
\]

\vspace{.2cm}

\begin{cor} For $p\in(0,\infty)$ the map $L:\Lp\rightarrow [-\infty,\infty[$ given by
  \begin{equation}
    L(T) = \tau(\log|T|), \qquad T\in\Lp,
  \end{equation}
is upper semicontinuous w.r.t. $\|\cdot\|_p$.
\end{cor}

\proof Indeed, this follows from Lemma~\ref{Lp-cont}, since for every
$T\in\Lp$ we have that
\[
L(T) = \inf_{\eps>0}L_\eps(T). \qquad \endproof
\]

\section{Unbounded $R$--diagonal operators}

Consider a von Neumann algebra $\CM$ equipped with a faithful, normal, tracial state $\tau$.

\begin{definition} For $T\in\tilde\CM$ with polar decomposition
$T=U|T|$, we denote by $W\cc(T)$ the von Neumann algebra generated by
$U$ and all the spectral projections of $|T|$.
\end{definition} 

Note that $T$ is affiliated with $W\cc(T)$ and that $W\cc(T)$ is the
smallest von Neumann subalgebra of $\CM$ with this property.

If $\CM_1$ and $\CM_2$ are finite von
  Neumann algebras with faithful, normal, tracial states $\tau_1$ and
  $\tau_2$, respectively, then any
  $\ast$--isomorphism $\phi:\CM_1\rightarrow \CM_2$ with
  $\tau_1=\tau_2\circ\phi$ is continuous w.r.t. the measure topologies
  on the two von Neumann algebras and thus has a unique
  extension to a (surjective) $\ast$--isomorphism
  $\tilde\phi:\tilde\CM_1\rightarrow\tilde\CM_2$.

\begin{definition}\label{def *-free} Let $S,T\in\tilde\CM$.
\begin{itemize}
  \item[(a)] We say that $S$ and $T$ have the same
  $\ast$--distribution, in symbols $S\underset{\ast\CD}{\sim} T$, if there exists a trace--preserving
  $\ast$--isomporphism $\phi$ from $W\cc(S)$ onto $W\cc(T)$ with
  $\tilde\phi(S)=T$.
  \item[(b)] We say that $S$ and $T$ are $\ast$--free if
  $W\cc(S)$ and $W\cc(T)$ are $\ast$-free.
\end{itemize}
\end{definition}

Note that in case $S$ and $T$ are bounded, the two definitions
  (a) and (b) given above agree with the ones given in \cite{VDN} . 

\vspace{.2cm}

Recall from \cite[p.~155~ff.]{NS} that if $U,H\in\CM$ are $\ast$-free
elements with $U$ Haar unitary, then $UH$ is $R$-diagonal in the sense of
Nica and Speicher (cf. \cite{NS}). Conversely, if $T\in \CM$ is
$R$--diagonal, then $T$ has the same $\ast$--distribution as a product
$UH$, where $U$ and $H$ are $\ast$--free elements in some tracial
$C\cc$--probability space, $U$ is a
Haar unitary, and $H\geq 0$. We therefore define
$R$--diagonality for operators in $\tilde\CM$ as follows: 

\begin{definition}\label{def R-diag} $T\in \tilde\CM$ is said to be
  \emph{$R$-diagonal} if there exist a  von Neumann algebra $\CN$, with a faithful,
  normal, tracial state, and $\ast$--free elements $U$ and $H$ in 
  $\tilde\CN$, such that $U$ is Haar unitary, $H\geq 0$, and such that
  $T$ has the same $\ast$--distribution as $UH$.
\end{definition}

\vspace{.2cm}

\begin{remark} Note that if $T\in\tilde\CM$ is $R$--diagonal with
  ${\rm ker}(T)=0$, then the partial isometry $V$ in the polar
  decomposition of $T$, $T=V|T|$, is a unitary ($\CM$ is finite). It
  follows from Definition~\ref{def R-diag} and Definition~\ref{def
  *-free} that $V$ is in fact a Haar
  unitary which is $\ast$--free from $|T|$.
\end{remark}

\vspace{.2cm}

In this section we will see that certain algebraic operations on (sets of
$\ast$-free) $R$--diagonal operators preserve $R$--diagonality, exactly as in
the bounded case (cf. \cite{HL}). Our proofs are to a large extent inspired
by the techniques used in \cite{HL} and in \cite{FL}. In particular,
we will repeatedly make use of \cite[Lemma~3.7]{HL} which we state here for the
convenience of the reader:

\begin{lemma}\cite{HL}\label{lemmaHL} Let $U\in\CM$ be a Haar unitary, and suppose $\CS\subset
  \CM$ is a set which is $\ast$-free from $U$. Then for any $n\in\N$,
  \begin{itemize}
    \item[(i)] the sets $\CS$, $U\CS U\cc$,.... are $\ast$-free,
    \item[(ii)] the sets $\CS$, $U\CS U\cc$,...., $U^{n-1}\CS (U\cc)^{n-1}$,
    $\{U^n\}$ are $\ast$-free,
    \item[(iii)] the sets $U\CS U\cc$,....,$U^{n}\CS (U\cc)^{n}$,
    $\{U^n\}$ are $\ast$-free.
  \end{itemize}
\end{lemma}

\vspace{.2cm}

\begin{prop}\label{Rdiag1} If $T\in \tilde\CM$ is $R$-diagonal with
    ${\rm ker}(T)=0$, then $T$ has an
  inverse $T^{-1}\in\tilde\CM$, and $T^{-1}$ is $R$-diagonal as
  well.
\end{prop}

\proof Let $T=V|T|$ be the polar decomposition of $T$ with $V\in\CM$
Haar unitary and $\ast$--free from $|T|$. Since ${\rm ker}(T)=0$, $T$
has an inverse $T^{-1}\in\tilde\CM$:
\[
T^{-1} = V\cc V |T|^{-1} V\cc = V\cc (V|T|V\cc)^{-1},
\]
where $V\cc$ is Haar unitary and, according to Lemma~\ref{lemmaHL}, it is
$\ast$--free from $V|T|V\cc$ and thus from $(V|T|V\cc)^{-1}$. This
shows that $T^{-1}$ is $R$--diagonal. $\endproof$

\vspace{.2cm}

\begin{lemma}\label{R-diag3}Let $S,T\in\tilde\CM$, and let $V\in\CM$
  be a Haar unitary. If $S,T$ and
  $V$ are $\ast$-free, then $VS$ and $TVS$ are $R$--diagonal.
\end{lemma}

\proof The case where $S$ and $T$ are bounded was treated by
F.~Larsen (cf. \cite[Lemma~3.6]{FL}). Our proof resembles the one given by
F.~Larsen.

Enlarging the algebra if necessary, we may assume that there are
Haar unitaries $V_1, V_2\in\CM$, such that $V_1, V_2$ and $S$ are
$\ast$-free and $V=V_1V_2$.

Since $W\cc(S)\subseteq \CM$ is finite, there is a unitary $U_1\in W\cc(S)$ such
that $S=U_1|S|$. Then $VS= V_1(V_2U_1)|S|$, where 
\begin{itemize}
  \item[(i)] $V_1$ is $\ast$-free from $|S|$ and $V_2U_1$,
  \item[(ii)] $\tau(V_1)=\tau(V_1\cc)=0$ and
  $\tau(V_2U_1)=\tau((V_2U_1)\cc)=0$,
  \item[(iii)] for all $A\in W\cc(|S|)$ with $\tau(A)=0$,
  $\tau(V_2U_1A)=\tau(V_2)\tau(U_1A)=0$, $\tau(AU_1\cc V_2\cc)=\tau(AU_1\cc)\tau(V_2\cc)=0$ and $\tau(V_2 U_1 A (V_2U_1)^{-1})=\tau(A)=0$.
\end{itemize}
It follows now from \cite[Lemma~2.4]{V1} that $V_1(V_2U_1)$ is $\ast$-free from
$|S|$. Thus, if $V_1(V_2U_1)$ is Haar unitary, then $S$ is $R$--diagonal.
Since $V_1$ is $\ast$-free
from $V_2U_1$, we get from \cite[Lemma~3.7]{HL} that for every $n\in\N$, the
operators
\[
V_1^n, V_1^{n-1}(V_2U_1)V_1^{1-n}, \ldots, V_1(V_2U_1)V_1^{-1}, V_2U_1
\]
are $\ast$-free. Consequently,
\begin{eqnarray*}
\tau((V_1V_2U_1)^n) &=& \tau(V_1^n[V_1^{1-n}(V_2U_1)V_1^{n-1}][V_1^{2-n}(V_2U_1)V_1^{n-2}]\cdots
[V_1^{-1}(V_2U_1)V_1]V_2U_1)\\
&=& \tau(V_1^n)\tau([V_1^{1-n}(V_2U_1)V_1^{n-1}][V_1^{2-n}(V_2U_1)V_1^{n-2}]\cdots
[V_1^{-1}(V_2U_1)V_1]V_2U_1) \\
&=&0.
\end{eqnarray*}
Then $\tau((V_1V_2U_1)^{-n})=\overline{\tau((V_1V_2U_1)^n)}=0$, and
$V_1V_2U_1$ is Haar unitary. Therefore $VS = V_1V_2U_1|S|$ is $R$--diagonal.

Now, $TVS =V(V\cc TVS)$. Put
\[
\CB_1 = W\cc(V),\;\; \CB_2=W\cc(T),\;\; {\rm and}\;\; \CB_3=W\cc(S).
\]
Then $\CB_1$,   $\CB_2$ and
$\CB_3$ are $\ast$-free. We may write $T$ as $T=U_2|T|$ for a unitary $U_2\in\CB_2$. Then
\begin{equation}
V\cc TV =(V\cc U_2 V)V\cc |T|V,
\end{equation}
where $V\cc|T|V$ is affiliated with $V\cc \CB_2 V$.

$\CB_3$ and $V\cc\CB_2 V$ are $\ast$-free, and according to \cite[Lemma~3.7]{HL}, $\CB_1$ and $V\cc\CB_2 V$ are
$\ast$-free. But then
$V$ is $\ast$-free from $\CB_4=\CB_3 \vee V\cc\CB_2 V$.

Since $S$ and $V\cc TV$ are both affiliated with $\CB_4$, their product,
$V\cc T VS$, is affiliated with $\CB_4$, so $V$ is $\ast$-free from
$V\cc TV S$. It follows now
from the first part of the proof that $TVS =V(V\cc TVS)$ is $R$-diagonal. $\endproof$

\vspace{.2cm}

\begin{prop}\label{Rdiag2} If $S,T\in\tilde\CM$ are $\ast$-free
  $R$-diagonal elements, then $ST$ is $R$-diagonal as well. Moreover,
\begin{equation}\label{free prod meas}
  \mu_{(ST)\cc ST}= \mu_{S\cc S}\boxtimes \mu_{T\cc T}.
\end{equation}
\end{prop}

\proof Taking a free product of tracial von Neumann algebras if
necessary, we can find a von Neumann algebra $\CN$ with faithful,
normal, tracial state $\omega$ and $\ast$--free elements
$U_1,H_1,U_2,H_2\in\tilde\CN$ such that $U_1,U_2$ are Haar unitaries,
$H_1,H_2\geq 0$, and $S\underset{\ast\CD}{\sim} U_1H_1$ and $T\underset{\ast\CD}{\sim} U_2H_2$.

Choose trace--preserving $\ast$--isomorphisms
\begin{eqnarray*}
  \phi_1: W\cc(S)&\rightarrow & W\cc(U_1H_1),\\
  \phi_2: W\cc(T)&\rightarrow & W\cc(U_2H_2),
\end{eqnarray*}
with $\tilde{\phi_1}(S)=U_1H_1$ and
$\tilde{\phi_2}(T)=U_2H_2$. $\phi_1$ and $\phi_2$ give rise to a
trace--preserving $\ast$--isomorphism 
\[
\phi=\phi_1\ast \phi_2: W\cc(S)\ast W\cc(T)\rightarrow
W\cc(U_1H_1)\ast W\cc(U_2H_2)
\]
(the free products are taken within the category of tracial von Neumann algebras)
with
\[
\tilde\phi(ST)= \tilde{\phi_1}(S)\tilde{\phi_2}(T)=U_1H_1U_2H_2.
\]
Thus, $\psi:=\phi|_{W\cc(ST)}$ is a trace--preserving
$\ast$--isomorphism onto $W\cc(U_1H_1U_2H_2)$ with
$\tilde\psi(ST)=U_1H_1U_2H_2$. According to Lemma~\ref{R-diag3},
$U_1(H_1U_2H_2)$ is $R$--diagonal, and hence $ST$ is $R$--diagonal. 

In order to prove \eqref{free prod meas}, note that if $S=0$, then
$\mu_{S\cc S}=\delta_0$, so that by the definition of multiplicative
free convolution given on p.~744 in \cite{BV},
\[
\mu_{S\cc S}\boxtimes \mu_{T\cc T} = \delta_{0}\boxtimes \mu_{T\cc T} =\delta_{0}.
\]
This shows that $\mu_{S\cc S}\boxtimes \mu_{T\cc T}= \mu_{(ST)\cc ST}$ if
$S=0$. The same holds if $T=0$. 

Now assume that $S,T\neq 0$. Note that
\begin{eqnarray*}
S\cc S& \underset{\ast\CD}{\sim} & H_1^2,\\
T\cc T& \underset{\ast\CD}{\sim} & H_2^2,\\
(ST)\cc ST & \underset{\ast\CD}{\sim}&  H_2U_2\cc
H_1^2 U_2H_2.
\end{eqnarray*}

Thus, \eqref{free prod meas} holds if
\[
\mu_{H_2U_2\cc
H_1^2 U_2H_2}= \mu_{H_1^2}\boxtimes \mu_{H_2^2}.
\]

For every $n\in\N$, the bounded
operators
\[
S_n=U_1\,H_1\, 1_{[0,n]}(H_1)
\]
and
\[
T_n=U_2\,H_2\, 1_{[0,n]}(H_2)
\]
are $\ast$--free. According to  \cite[Lemma~3.9]{HL} they
are both $R$-diagonal in the sense of Nica and
Speicher (cf. \cite{NS}). Then, by \cite[Proposition~3.6]{HL},
\begin{eqnarray}
  \mu_{(S_n T_n)\cc S_n T_n}&=& \mu_{S_n\cc S_n}\boxtimes \mu_{T_n\cc T_n}\label{eq6-10}.
\end{eqnarray}

Since $S_n\rightarrow U_1H_1$ and $T_n\rightarrow U_2H_2$ in the
measure topology,
$(S_nT_n)\cc S_nT_n \rightarrow H_2U_2\cc
H_1^2 U_2H_2$
in measure as well. These facts imply that $ \mu_{S_n\cc S_n}\overset{w\cc}{\rightarrow}
\mu_{H_1^2} $, $ \mu_{T_n\cc T_n}\overset{w\cc}{\rightarrow}
\mu_{H_2^2} $ and $\mu_{(S_nT_n)\cc S_nT_n}\overset{w\cc}{\rightarrow}
\mu_{H_2U_2\cc H_1^2 U_2H_2} $  . Moreover,  $\mu_{H_1^2}\neq
\delta_0$ and  $\mu_{H_2^2}\neq
\delta_0$, because $S\cc S$ and $T\cc T$ are non--zero. Hence, by
\cite[Corollary~6.7]{BV} and by \eqref{eq6-10},
\begin{equation*}
\mu_{H_2U_2\cc H_1^2 U_2H_2} = w\cc -\lim_{n\rightarrow\infty}\mu_{S_n\cc S_n}\boxtimes
\mu_{T_n\cc T_n} =  \mu_{H_1^2}\boxtimes \mu_{H_2^2}. \endproof
\end{equation*}

\vspace{.2cm}

\begin{prop}\label{Rdiag4} Let $S\in\tilde\CM$ be $R$-diagonal, and let
  $n\in\N$. Then $S^n$ is $R$-diagonal. Moreover,
  \begin{equation}\label{30april}
  \mu_{(S^n)\cc S^n}= \mu_{S\cc S}^{\boxtimes n}.
  \end{equation}
\end{prop}

\proof Choose a von Neumann algebra $\CN$ with faithful, normal,
tracial state $\omega$ and with $\ast$--free elements
$U,H\in\tilde\CN$ such that $U$ is Haar unitary, $H\geq 0$, and
$S\underset{\ast\CD}{\sim} UH$. Then $S^n \underset{\ast\CD}{\sim} (UH)^n$. Since
\[
(UH)^n = U^n[U^{1-n}HU^{n-1}] [U^{2-n}HU^{n-2}]\cdots [U^{-1}HU]H,
\]
where
\[
U^n, \; U^{1-n}HU^{n-1}, \; U^{2-n}HU^{n-2},\ldots,\; U^{-1}HU,\; H
\]
are $\ast$--free (cf. Lemma~\ref{lemmaHL}~(ii)), and $U^n$ is Haar
unitary, Lemma~\ref{R-diag3} gives us that $(UH)^n$ is $R$--diagonal, and
hence $S^n$ is.

In order to prove \eqref{30april}, note that if $\mu_{S\cc
  S}=\delta_0$, then $S=S^n=0$ and \eqref{30april} trivially
  holds.

Now assume that $\mu_{S\cc S}\neq \delta_0$. For $k\in\N$ define $S_k\in\CM$ and $T_k\in \CN$ by
\[
S_k=S\,1_{[0,k]}(|S|)\qquad {\rm and} \qquad T_k = U\,H\,1_{[0,k]}(H).
\]
Then  $T_k \underset{\ast\CD}{\sim} S_k$. Moreover, by Lemma~\ref{R-diag3}, $T_k$ is $R$-diagonal in the sense
of Nica and Speicher, so $S_k$ is $R$--diagonal. It now follows from
\cite[Proposition~3.10]{HL} that
\begin{equation}\label{eq6-12}
\mu_{[(S_k)^n]\cc(S_k)^n}=\mu_{[(T_k)^n]\cc(T_k)^n} =\mu_{T_k\cc T_k}^{\boxtimes n}=\mu_{S_k\cc S_k}^{\boxtimes n}.
\end{equation}

As $k$ tends to infinity, $S_k\cc
S_k \rightarrow S\cc S$ and $[(S_k)^n]\cc(S_k)^n \rightarrow (S^n)\cc S^n$
in the measure topology. Since $\mu_{S\cc S}\neq\delta_0$, we infer from
\cite[Corollary~6.7]{BV} and from \eqref{eq6-12} that
\[
\mu_{(S^n)\cc S^n} =w\cc-\lim_{k\rightarrow\infty}
  \mu_{[(S_k)^n]\cc(S_k)^n} = w\cc-\lim_{k\rightarrow\infty} \mu_{S_k\cc
  S_k}^{\boxtimes n} =\mu_{S\cc S}^{\boxtimes n}. \endproof
\]

\vspace{.2cm}

\begin{definition} For $\mu\in\Prob(\R, \B)$ let $\tilde\mu$ denote the
  {\it symmetrization} of $\mu$. That is, $\tilde\mu\in \Prob(\R, \B)$ is given
  by
  \[
  \tilde\mu(B) = \textstyle{\frac12}(\mu(B)+\mu(-B)), \qquad (B\in\B).
  \]
\end{definition}

\vspace{.2cm}

\begin{prop}\label{Rdiag5} Let $S,T\in\tilde\CM$ be $\ast$--free $R$--diagonal
  elements. Then
  \begin{equation}
    \tilde\mu_{|S+T|}= \tilde\mu_{|S|}\boxplus \tilde\mu_{|T|}.
  \end{equation}
\end{prop}

\proof As in the proof of Proposition~\ref{Rdiag2}, choose
$(\CN,\omega)$ and $\ast$--free elements $U_1,H_1,U_2,H_2\in\tilde\CN$ such that $U_1,U_2$ are Haar unitaries,
$H_1,H_2\geq 0$, and $S\underset{\ast\CD}{\sim} U_1H_1$ and $T\underset{\ast\CD}{\sim} U_2H_2$.

Again, for $n\in\N$, let
\[
S_n=U_1\,H_1\, 1_{[0,n]}(H_1)
\]
and
\[
T_n=U_2\,H_2\, 1_{[0,n]}(H_2).
\]
Then $S_n$ and $T_n$ are $\ast$--free and $R$--diagonal
and therefore, according to \cite[Proposition~3.5]{HL},
\begin{equation}
  \tilde\mu_{|S_n+T_n|} =  \tilde\mu_{|S_n|}\boxplus \tilde\mu_{|T_n|}.
\end{equation}

$|S_n|\rightarrow H_1$ and $|T_n|\rightarrow H_2$ in measure, implying that
$\mu_{|S_n|}\overset{w\cc}{\rightarrow} \mu_{H_1}=\mu_{|S|}$ and $\mu_{|T_n|}\overset{w\cc}{\rightarrow} \mu_{H_2}=\mu_{|T|}$
 . Then we also have weak convergence of the symmetrized measures:
\[
\tilde\mu_{|S_n|}\overset{w\cc}{\rightarrow} \tilde\mu_{|S|} \quad {\rm and}
\quad \tilde\mu_{|T_n|}\overset{w\cc}{\rightarrow} \tilde\mu_{|T|}.
\]

Let $d$ denote the L\'evy metric on $\Prob(\R, \B)$
(cf. \cite[p. 743]{BV}). Then $d$ induces the topology of weak convergence,
and according to \cite[Proposition~4.13]{BV} and the above observations,
\[
d(\tilde\mu_{|S|}\boxplus \tilde\mu_{|T|},\tilde\mu_{|S_n|}\boxplus
\tilde\mu_{|T_n|} )\leq d(\tilde\mu_{|S|}, \tilde\mu_{|S_n|})+
d(\tilde\mu_{|T|}, \tilde\mu_{|T_n|}) \rightarrow 0 \quad {\rm as}\;
n\rightarrow \infty.
\]
It follows that
\begin{eqnarray}
\tilde\mu_{|S|}\boxplus \tilde\mu_{|T|}  & = & w\cc -\lim_{n\rightarrow
  \infty}\tilde\mu_{|S_n|}\boxplus\tilde\mu_{|T_n|}\nonumber\\
& = & w\cc -\lim_{n\rightarrow\infty} \tilde\mu_{|S_n + T_n|}.
\end{eqnarray}

Since $S$ and $T$ ($U_1H_1$ and $U_2H_2$, resp.) are $\ast$--free with
$S\underset{\ast\CD}{\sim} U_1H_1$ and $T\underset{\ast\CD}{\sim}
U_2H_2$, it follows that $S+T\underset{\ast\CD}{\sim} U_1H_1+ U_2H_2$. Moreover,
$|S_n+T_n|\rightarrow |U_1H_1+U_2H_2|\underset{\ast\CD}{\sim}|S+T|$ in measure, and thus
$\tilde\mu_{|S_n + T_n|}\overset{w\cc}{\rightarrow} \tilde\mu_{|S + T|}$  . Finally,
this implies that
\[
\tilde\mu_{|S|}\boxplus \tilde\mu_{|T|} =  \tilde\mu_{|S + T|}. 
\endproof
\]


\vspace{.2cm}

We close this section by proving two simple results on the
$S$--transform of probability measures on $(0,\infty)$
(cf. \cite{BV}).

For $\mu\in\Prob((0,\infty), \B)$ define $\psi_\mu: \C\setminus
(0,\infty)\rightarrow \C$
by
\begin{equation}
  \psi_\mu(z)=\int_0^\infty \frac{1}{1-zt}\,\d\mu(t) - 1, \qquad z\in
  \C\setminus (0,\infty).
\end{equation}

Then $\psi_\mu$ is analytic and satisfies
\begin{itemize}
  \item[(i)] $\psi_\mu'(t)>0$, $t\in (-\infty, 0)$,
  \item[(ii)]$\psi_\mu(z)\rightarrow -1$ as $z\rightarrow  -\infty$,
  \item[(iii)]$\psi_\mu(z)\rightarrow 0$ as $z\rightarrow 0$.
\end{itemize}

Hence, $\psi_\mu$ maps a (connected) neighbourhood $\CU_\mu$ of $(-\infty, 0)$ injectively
onto a neighbourhood $\CV_\mu$ of $(-1, 0)$. Define $\chi_\mu$, $\CS_\mu:
\CV_\mu\rightarrow \C$ by
\begin{eqnarray}
\chi_\mu(z) & = & \psi_\mu^{-1}(z), \qquad z\in \CV_\mu,\\
\CS_\mu(z)&=&\frac{z+1}{z}\chi_\mu(z), \qquad z\in \CV_\mu.
\end{eqnarray}

\vspace{.2cm}

\begin{prop}\label{S1} The map $\mu\mapsto \CS_\mu$ is one-to-one on
  $\Prob((0,\infty), \B)$.
\end{prop}

\proof Suppose $\mu, \nu\in \Prob((0,\infty),\B)$ with $\CS_\mu =
\CS_\nu$. That is, in a neighbourhood $\CV=\CV_\mu\cap \CV_\nu $ of  $(-1,
0)$, $\chi_\mu$ agrees with $\chi_\nu$. It follows that on $(-\infty,0)$,
$\psi_\mu$ agrees with $\psi_\nu$, and then, by uniqueness of analytic continuation,
\begin{equation}
 \psi_\mu(\textstyle{\frac1\lambda})=\psi_\nu(\textstyle{\frac1\lambda}), \qquad
 \lambda\in\C \setminus [0,\infty[.
\end{equation}
That is, the Stieltjes-transforms $G_\mu$ and $G_\nu$ agree on $\C
\setminus [0,\infty[ $. Recall that
\begin{equation}
\d\mu(x)=-\textstyle{\frac1\pi}\lim_{y\rightarrow 0^+}G_\mu(x+\i y)\d x
\end{equation}
(weak convergence of measures), and similarly,
\begin{equation}
\d\nu(x)=-\textstyle{\frac1\pi}\lim_{y\rightarrow 0^+}G_\nu(x+\i y)\d x.
\end{equation}
Thus $\mu = \nu$. $\endproof$

\vspace{.2cm}

\begin{prop}\label{S4} Let $\CM$ be a II$_1$-factor with tracial state
  $\tau$, and let $a\in \tilde\CM_+$ with ${\rm
    ker}(a)=\{0\}$. Then for all $z$ in a neighbourhood of $(-1,0)$,
  \begin{equation}
    \CS_{\mu_{a^{-1}}}(z)= \frac{1}{\CS_{\mu_{a}}(-1-z)}.
  \end{equation}
\end{prop}

\proof Let $z\in \C\setminus [0,\infty[$. Then
\begin{eqnarray*}
  \psi_{a^{-1}}(z)&=& \int_0^\infty \frac{1}{1-zt}\,\d\mu_{a^{-1}}(t) -1\\
  & = &  \int_0^\infty \frac{1}{1-\frac zt}\,\d\mu_{a}(t) -1\\
  & = & \int_0^\infty \frac{z}{t-z}\d\mu_a(t),
\end{eqnarray*}
and hence
\begin{equation}
\psi_{a^{-1}}(\textstyle{\frac1z})= -\int_0^\infty
\frac{1}{1-zt}\,\d\mu_a(t) = -(\psi_a(z)+1).
\end{equation}
It follows that for all $z\in \C\setminus
 [0,\infty[$,
\begin{equation}
z= \chi_a(\psi_a(z))= \chi_a(-1-\psi_{a^{-1}}(\textstyle{\frac1z})),
\end{equation}
implying that $w = \psi_{a^{-1}}(\textstyle{\frac1z})$ satisfies
\begin{equation}
\chi_{a^{-1}}(w)= \frac1z = \frac{1}{\chi_a(-1-w)},
\end{equation}
and thus
\begin{equation}\label{eq6-1}
\CS_{\mu_{a^{-1}}}(w)\cdot\CS_{\mu_{a}}(-1-w)=1.
\end{equation}
\eqref{eq6-1} holds for all $w\in \psi_{a^{-1}}(\C\setminus [0,\infty[)$
and in particular for all $w$ in a neighbourhood of $(-1,0)$. $\endproof$

\section{The Brown measure of an unbounded $R$--diagonal operator}

The Brown measure of a general bounded $R$--diagonal operator was
computed in \cite[Theorem~4.4]{HL}. We will genralize this result to
unbounded $R$--diagonal elements in $\CMD$. Our proof will take a
different route than the one in \cite{HL}. This new approach will
enable us to obtain an estimate of the $p$--norm of the
resolvent $(T-\lambda\unit)^{-1}$, $0<p<1$, for special $R$--diagonal elements $T$ (cf. Section~5). 

\begin{lemma}\label{Bmeas1} Let $T\in\tilde\CM$ be an $R$--diagonal element, and let
  $U\in\CM$ be a Haar unitary which is $\ast$--free from $T$. Then for
  every $\lambda\in\C$,
  \begin{equation}
    |T-\lambda\unit| \underset{\ast\CD}{\sim} |T+|\lambda|U|.
  \end{equation}
\end{lemma}

\proof By passing to a larger algebra, we may assume that $T=V|T|$
where $V\in\CM$ is a Haar unitary and $U,V$ and $|T|$ are
$\ast$--free. The case $\lambda=0$ is trivial. For $\lambda\neq 0$,
let $\alpha= -\frac{\lambda}{|\lambda|}$. Then $\alpha U\cc V$ is a
Haar unitary which is $\ast$--free from $T$. Hence,
\[
\alpha U\cc V|T| \underset{\ast\CD}{\sim} T.
\]
Therefore,
\begin{eqnarray*}
  |T-\lambda\unit| & \underset{\ast\CD}{\sim} & |\alpha U\cc
   V|T|-\lambda\unit|\\
   &=& |T-\overline\alpha \lambda U|\\
   &=& |T+|\lambda|U|.
\end{eqnarray*}

\vspace{.2cm}

\begin{lemma}\label{lemma h} Let $T\in\tilde\CM$ be an $R$--diagonal operator, and
  define
  \[
  h(s)= s\,\tau\big((T\cc T+s^2\unit)^{-1}\big), \qquad s>0.
  \]
  Moreover, for $\lambda\in\C\setminus\{0\}$, set
  \[
  h_\lambda(s) = s\,\tau\big([(T-\lambda\unit)\cc
  (T-\lambda\unit)+s^2\unit]^{-1}\big).
  \]
  Then there exists an $s_\lambda>0$ such that for $s>s_\lambda$,
  \[
  h(s)=
  h_\lambda\Bigg(s+\frac{\sqrt{1-4|\lambda|^2h(s)^2}-1}{2h(s)}\Bigg).
  \]
\end{lemma}

\proof By passing to a larger algebra, we may assume that there exists
a Haar unitary $U\in\CM$ which is $\ast$--free from $T$. Then,
according to Lemma~\ref{Bmeas1},
\[
|T-\lambda\unit| \underset{\ast\CD}{\sim}|T+|\lambda|U|.
\]
It follows now from Proposition~\ref{Rdiag5} that
\[
\tilde{\mu}_{|T-\lambda\unit|} = \tilde{\mu}_{|T|}\boxplus
\tilde{\mu}_{|\lambda|\unit}= \tilde{\mu}_{|T|}\boxplus \nu,
\]
where $\nu= \frac12 (\delta_{-|\lambda|}+\delta_{|\lambda|})$.

For $\beta >0$ define
\[
\Omega_\beta = \{w\in \C \,|\, \, 0<|w|<\beta, \,
\textstyle{\frac{5\pi}{4}} < \arg (w) < \textstyle{\frac{7\pi}{4}}\}.
\]
According to \cite[Corollary~5.8]{BV}, there is a $\beta>0$ such that for
every $w\in \Omega_\beta$,
\[
\CR_{\tilde\mu_{|T-\lambda\unit|}}(w)= \CR_{\tilde\mu_{|T|}}(w) + \CR_\nu(w),
\]
where
\[
\CR_\nu(w)= \frac{\sqrt{1+4|\lambda|^2w^2}-1}{2w},
\]
and
\[
G_{\tilde\mu_{|T|}}(\i s)= -\i h(s), \qquad s>0,
\]
whence
\[
\CR_{\tilde\mu_{|T|}}(-\i h(s)) + \frac{1}{-\i h(s)}=
G_{\tilde\mu_{|T|}}^{<-1>}(-\i h(s))=\i s,  \qquad s>0.
\]

Take $s_\lambda>0$ such that for every $s>s_\lambda$, $-\i h(s)\in \Omega_\beta$. Then,
when $s>s_\lambda$,
\[
\CR_{\tilde\mu_{|T-\lambda\unit|}}(-\i h(s)) = \i s +  \frac{1}{\i h(s)} +
\frac{\sqrt{1-4 |\lambda|^2h(s)^2}-1}{-2\i h(s)},
\]
implying that
\[
h(s)= h_\lambda \Bigg(s+\frac{\sqrt{1-4|\lambda|^2h(s)^2}-1}{2h(s)}\Bigg).
\]
That is, when $s>s_\lambda$ and
\[
t = s+\frac{\sqrt{1-4|\lambda|^2h(s)^2}-1}{2h(s)},
\]
then $h(s)=h_\lambda(t)$. $\endproof$

\vspace{.2cm}

Note that if
\[
t= s + \frac{\sqrt{1-4|\lambda|^2h(s)^2}-1}{2h(s)},
\]
then $(s,t)$ satisfies the following equation:
\begin{equation}
(s-t)\Big(\frac{1}{h(s)}-s+t\Big)=|\lambda|^2.
\end{equation}
In the following we will investigate this equation further.

\vspace{.2cm}

\begin{definition} Let $m,n\in\N$, and let $U$ be an open set in
  $\R^m$. A map $f:U\rightarrow \R^n$ is said to be \emph{analytic} if
  it has a power series expansion in $m$ variables in a neighborhood
  of every $x\in U$.
\end{definition}

\vspace{.2cm}

We shall need the following two well--known lemmas about analytic
functions of several variables:

\begin{lemma}\label{lemma5.1} Let $U$ be a connected, open subset of $\R^m$. If $f,g:
  U\rightarrow \R^n$ are two analytic functions which coincide on a
  non--empty, open subset $V$ of $U$, then $f=g$.
\end{lemma}

\vspace{.2cm}

\begin{lemma}\label{lemma5.3} Let $U\subseteq \R^m$ be open and let $f: U\rightarrow
  \R^m$ be an analytic function for which the Jacobian $\CJ(x_0)={\rm
  det} f'(x_0)$ is non--zero for some $x_0\in U$. Then $f$ is
  one--to--one in some neighborhood $V$ of $x_0$, and the inverse of
  $f|_V$ is analytic in a neighborhood of $f(x_0)$.
\end{lemma}

\vspace{.2cm}

\begin{lemma}\label{lemma5.4} Let $\mu$ be a probability measure on $[0,\infty)$, and
  define
  \begin{equation}
  h(s)=\int_0^\infty\frac{s}{s^2+u^2}\,\d\mu(u), \qquad s>0.
  \end{equation}
Then $h$ is analytic on $(0,\infty)$. Moreover, if $\mu$ is not a
Dirac measure, then for all $s>0$,
\[
0<h(s)<\frac1s \qquad {\rm and} \qquad h'(s)<\frac{h(s)}{s}-2h(s)^2.
\]
\end{lemma}  

\proof Since
\[
h(s)=\frac12 \int_0^\infty \Bigg(\frac{1}{s+\i u}+\frac{1}{s-\i
  u}\Bigg)\,\d\mu, \qquad s>0,
\]
$h$ has a complex analytic extension
\[
\tilde{h} : \{z\in\C\,|\, \im z>0\}\rightarrow \C
\]
given by the same formula. In particular, $h$ is an analytic function
of $s\in (0,\infty)$. If $\mu$ is not a Dirac measure, then $\mu\neq
\delta_0$, and so
$h(s)>0$ for all $s>0$. Moreover,
\[
s\,h(s) = \int_0^\infty\frac{s^2}{s^2+u^2}\,\d\mu(u)<1, \qquad s>0.
\]
Finally, for $s>0$,
\begin{eqnarray*}
  h(s)^2 &=& \int_0^\infty\int_0^\infty
  \frac{s}{s^2+u^2}\frac{s}{s^2+v^2}\,\d\mu(u)\,\d\mu(v)\\
  &\leq &  \int_0^\infty\int_0^\infty \frac12
  \Bigg(\Big(\frac{s}{s^2+u^2}\Big)^2 +
  \Big(\frac{s}{s^2+v^2}\Big)^2\Bigg)\,\d\mu(u)\,\d\mu(v)\\
  &=& \int_0^\infty \frac{s^2}{(s^2+u^2)^2}\,\d\mu(u)\\
  &=& \frac12 \Bigg(\int_0^\infty
  \frac{s^2+u^2}{(s^2+u^2)^2}\,\d\mu(u) + \int_0^\infty
  \frac{s^2-u^2}{(s^2+u^2)^2}\,\d\mu(u) \Bigg)\\
  &=& \frac12 \Bigg(\frac{h(s)}{s} -h'(s)\Bigg).
\end{eqnarray*}
Hence,
\[
h'(s)\leq \frac{h(s)}{s}-2 h(s)^2,
\]
and equality holds if and only if the product measure
$\mu\otimes\mu$ is concentrated on the diagonal $\{(u,u)\,|\,
u>0\}$. But this would imply that $\mu$ is a Dirac measure. Thus, if $\mu$ is not a
Dirac measure, then
\[
h'(s)< \frac{h(s)}{s}-2 h(s)^2, \qquad s>0\,\,\, \endproof
\]

\vspace{.2cm}

\begin{lemma}\label{lemma5.5} Let $\mu$ be a probability measure on $[0,\infty)$ which
  is not a Dirac measure, and
  put
  \[
  \lambda_1(\mu) = \Bigg(\int_0^\infty
  \frac{1}{u^2}\,\d\mu(u)\Bigg)^{-\frac12} \quad {\rm and} \quad \lambda_2(\mu) = \Bigg(\int_0^\infty
  u^2\,\d\mu(u)\Bigg)^{\frac12},
  \]
  with the convention that $\infty^{-\frac12}=0$. 
Then $0\leq \lambda_1(\mu)<\lambda_2(\mu)\leq \infty$. 
\end{lemma}

\proof Clearly, $\lambda_1(\mu)<\infty$, and since $\mu\neq \delta_0$,
$\lambda_2(\mu)>0$. The lemma is then trivially true if $\lambda_1(\mu)=0$
or $\lambda_2(\mu)= +\infty$. Thus, we can assume that
$\lambda_1(\mu), \lambda_2(\mu)\in (0,\infty)$. Then, by the Schwartz
inequality,
\begin{eqnarray*}
\frac{\lambda_2(\mu)}{\lambda_1(\mu)}&=& \Bigg(\int_0^\infty
u^2\,\d\mu(u)\Bigg)^\frac12 \Bigg(\int_0^\infty
\frac{1}{u^2}\,\d\mu(u)\Bigg)^\frac12\\
&\geq & \int_0^\infty u\,\frac1u\,\d\mu(u)\\
&=& 1,
\end{eqnarray*}
and equality holds if and only if for some $c\in (0,\infty)$, $\frac1u
= cu$ holds for $\mu$--a.e. $u\in [0,\infty)$. However, this can not
be the case when $\mu$ is not a Dirac measure. $\endproof$

\vspace{.2cm}

\begin{lemma}\label{lemma5.6} Let $\mu, \lambda_1(u)$ and $\lambda_2(\mu)$ be as in
  Lemma~\ref{lemma5.5}, and let $h$ be as in
  Lemma~\ref{lemma5.4}. Then put
  \[
  k(s,t)=(s-t)\Bigg(\frac{1}{h(s)}-s+t\Bigg), \qquad s>0,\, t\in\R.
  \]
  Then $k$ is an analytic function on $(0,\infty)\times\R$. Moreover,
  for $t>0$ the map $s\mapsto k(s,t)$ is a strictly increasing
  bijection of $(t,\infty)$ onto $(0,\infty)$, and for $t=0$ the map
  $s\mapsto k(s,t)$ is a strictly increasing bijection of $(0,\infty)$
  onto $(\lambda_1(\mu)^2,\lambda_2(\mu)^2)$.
  \end{lemma}

  \proof Clearly, $k$ is analytic. Moreover,
  \begin{equation}\label{5.2}
  \frac{\partial k}{\partial s}(s,t) = \frac{1}{h(s)}
  -(s-t)\Bigg(2+\frac{h'(s)}{h(s)^2}\Bigg).
  \end{equation}
For $s\in (0,\infty)$, we get from Lemma~\ref{lemma5.4} that
\[
\frac{\partial k}{\partial s}(s,0)=
\frac{s}{h(s)^2}\Bigg(\frac{h(s)}{s} -2h(s)^2 - h'(s)\Bigg)>0,
\]
and
\[
\frac{\partial k}{\partial s}(s,s)=\frac{1}{h(s)}>s.
\]
Since the right--hand side of \eqref{5.2} is an affine function of
$t\in\R$, it follows that
\begin{equation}\label{5.3}
  \frac{\partial k}{\partial s}(s,t)>t, \quad s>0, \quad t\in [0,s].
\end{equation}
Hence, $s\mapsto k(s,t)$ is a strictly increasing function of $s\in
(t,\infty)$ for every $t\in [0,\infty)$. For $s>t>0$,
\begin{equation}\label{5.6a}
k(s,t)= \int_t^s \frac{\partial k}{\partial s'}(s',t)\,\d s' >\int_t^s
t\,\d s' =t(s-t).
\end{equation}
Hence, when $t>0$,
\[
\lim_{s\rightarrow \infty}k(s,t)= \infty,
\]
and
\[
\lim_{s\rightarrow t+} k(s,t)=k(t,t)=0.
\]
Thus, $s\mapsto k(s,t)$ is a bijection of $(t,\infty)$ onto
$(0,\infty)$.

Next, consider the case $t=0$. We have already seen that $s\mapsto
k(s,0)$ is strictly increasing on $(0,\infty)$. Note that for $s>0$,
\[
k(s,0) = \frac{1-sh(s)}{h(s)/s} = \frac{n(s)}{d(s)}
\]
where
\[
n(s)= \int_0^\infty \frac{u^2}{s^2+u^2}\,\d\mu(u) \quad {\rm and}
\quad d(s)=\int_0^\infty \frac{1}{s^2+u^2}\,\d\mu(u).
\]
By the monotone convergence theorem,
\[
\lim_{s\rightarrow 0+}n(s)=1,
\]
\[
\lim_{s\rightarrow 0+} d(s) = \int_0^\infty \frac{1}{u^2}\,\d\mu(u) = \frac{1}{\lambda_1(\mu)^2},
\]
\[
\lim_{s\rightarrow \infty} s^2n(s) = \int_0^\infty u^2\,\d\mu(u) =
\lambda_2(\mu)^2,
\]
and
\[
\lim_{s\rightarrow \infty} s^2 d(s)=1.
\]
Hence,
\[
\lim_{s\rightarrow 0+} k(s,0)=\lambda_1(\mu)^2,
\]
and
\[
\lim_{s\rightarrow \infty}k(s,0)=\lambda_2(\mu)^2.
\]
This shows that $s\mapsto k(s,0)$ is a bijection of $(0,\infty)$ onto
$(\lambda_1(\mu)^2, \lambda_2(\mu)^2)$. $\endproof$

\vspace{.2cm}

\begin{definition}\label{def5.10} Let $\mu, \lambda_1(u)$ and $\lambda_2(\mu)$ be as in
  Lemma~\ref{lemma5.5}, let $h$ be as in
  Lemma~\ref{lemma5.4}, and let $k$ be as in Lemma~\ref{lemma5.6}. For $\lambda, t\in (0,\infty)$, let $s(\lambda, t)$ denote the unique solution $s\in (t,\infty)$ to the equation $k(s,t)=\lambda^2$ (cf. Lemma~\ref{lemma5.6}), and for $\lambda\in (\lambda_1(\mu), \lambda_2(\mu))$, let $s(\lambda,0)$ denote the unique solution $s\in (0,\infty)$ to the equation $k(s,0)=\lambda^2$.
  \end{definition}

\vspace{.2cm}

\begin{lemma}\label{lemma5.8} The function $(\lambda,t)\mapsto s(\lambda,t)$ is analytic in $(0,\infty)\times (0,\infty)$. Moreover, for $\lambda\in (\lambda_1(\mu),\lambda_2(\mu))$, 
\begin{equation}\label{5.4}
\lim_{t\rightarrow 0+}s(\lambda,t)=s(\lambda,0).
\end{equation}
\end{lemma}

\proof Let
\[
\Omega = \{(s,t)\in\R^2\,|\, 0<t<s\}.
\]
According to Lemma~\ref{lemma5.6}, $k$ is a strictly positive, analytic function in $\Omega$. Let
\[
F(s,t)=(\sqrt{k(s,t)}, t), \quad (s,t)\in\Omega.
\]
Then $F$ is analytic in $\Omega$, and by Lemma~\ref{lemma5.6}, $F$ is a one--to--one map of $\Omega$ onto $(0,\infty)\times (0,\infty)$. Moreover, its inverse $F^{-1}: (0,\infty)\times (0,\infty)\rightarrow \Omega$ is given by
\[
F^{-1}(\lambda,t)= (s(\lambda,t),t), \quad s,t>0.
\]
The Jacobian of $F$ is 
\[
\CJ(F)(s,t)= \frac{\partial}{\partial s}\sqrt{k(s,t)} = \frac{1}{2\sqrt{k(s,t)}}\frac{\partial k}{\partial s}(s,t),
\]
which by \eqref{5.3} is strictly positive for all $(s,t)\in\Omega$. Hence, by Lemma~\ref{lemma5.3}, $F^{-1}$ is analytic in $(0,\infty)\times (0,\infty)$. In particular, $s(\lambda,t)$ is analytic in $(0,\infty)\times (0,\infty)$. 

Now, let $\lambda_0\in(\lambda_1(\mu),\lambda_2(\mu))$ and put $s_0=s(\lambda_0,0)$. Then $k(s_0,0)=\lambda_0^2$, and by the proof of Lemma~\ref{lemma5.6}, $\frac{\partial k}{\partial s}(s_0,0)>0$. Let
\[
F_0(s,t):= (\sqrt{k(s,t)},t).
\]
$F_0$ is then analytic in some neighborhood $U_0$ of
$(s_0,0)$. Moreover, $\CJ(F_0)(s_0,0)\neq 0$, and therefore, by
Lemma~\ref{lemma5.3}, $F_0$ has an analytic inverse $F_0^{-1}$ in a
neighborhood $V_0$ of $F_0(s_0,0)=(\lambda_0,0)$. Clearly,
$F_0^{-1}(\lambda,t)=F^{-1}(\lambda,t)$, whenever $(\lambda,t)\in
V_0\cap [(0,\infty)\times (0,\infty)]$, and then and $F_0^{-1}(\lambda,t)\in\Omega$. 

Note that 
\begin{equation}\label{5.5}
\lim_{t\rightarrow 0+}F_0^{-1}(\lambda_0,t)= F_0^{-1}(\lambda_0,0)=(s_0,0), 
\end{equation}
and since the second coordinate of $F_0^{-1}(\lambda_0,t)$ is $t$, we conclude that $F_0^{-1}(\lambda_0,t)\in\Omega$, eventually as $t\rightarrow 0+$. Hence,
\begin{eqnarray*}
(s_0,0)&=& \lim_{t\rightarrow 0+}F_0^{-1}(\lambda_0,t)\\
&=& \lim_{t\rightarrow 0+}F^{-1}(\lambda_0,t)\\
&=& \lim_{t\rightarrow 0+}(s(\lambda_0,t),t),
\end{eqnarray*}
and therefore,
\[
\lim_{t\rightarrow 0+}s(\lambda_0,t)=s_0=s(\lambda_0,0). \qquad \endproof
\]

\vspace{.2cm}

\begin{remark}\label{rem5.11} We get from Lemma~\ref{lemma5.8} that
  \begin{equation}\label{5.6}
  \lim_{t\rightarrow 0+}s(\lambda,t)=0, \qquad 0<\lambda\leq\lambda_1(\mu),
\end{equation}
and
\begin{equation}\label{5.7}
  \lim_{t\rightarrow 0+}s(\lambda,t)=+\infty, \qquad
  \lambda\geq\lambda_2(\mu). 
\end{equation}
Indeed, for fixed $t>0$, $\lambda\mapsto s(\lambda,t)$ is a monotonically
increasing function of $\lambda$. Hence, if
$0<\lambda\leq \lambda_1(\mu)$, then
\[
\limsup_{t\rightarrow 0+}s(\lambda,t)\leq \limsup_{t\rightarrow
  0+}s(\lambda',t)=s(\lambda',0),
\]
for all $\lambda'\in (\lambda_1(\mu),\lambda_2(\mu))$.

But $\lambda' \mapsto s(\lambda',0)$ is the inverse function of
$s\mapsto \sqrt{k(s,0)}$, and hence $\lambda'\mapsto s(\lambda',0)$ is
a bijection of $(\lambda_1(\mu),\lambda_2(\mu))$ onto $(0,\infty)$. It
follows that $\limsup_{t\rightarrow 0+}
s(\lambda,t)=0$, and this proves \eqref{5.6}.

For
$\lambda\geq \lambda_2(\mu)$, a similar argument shows that $\liminf_{t\rightarrow 0+}
s(\lambda,t)=+\infty$, and this proves \eqref{5.7}.
\end{remark}

\vspace{.2cm}

\begin{lemma}\label{lemma5.12}
  Let $\lambda>0$. Then
  \begin{itemize}
    \item[(i)] $\lim_{t\rightarrow\infty}(s(\lambda,t)-t)=0$, and
    \item[(ii)] there exists a $t_\lambda >0$ such that when
    $t>t_\lambda$ and $s=s(\lambda,t)$, then
    \[
    t=s+\frac{\sqrt{1+4\lambda^2h(s)^2}-1}{2h(s)}.
    \]
  \end{itemize}
\end{lemma}

\proof Fix $t>0$, and put $s=s(\lambda,t)$. Then by
Definition~\ref{def5.10}, $s>t$ and $k(s,t)=\lambda^2$. According to
\eqref{5.6a}, $k(s,t)>t(s-t)$. Hence,
\[
0<s-t<\frac{\lambda^2}{t}.
\]
This proves (i). With $s$ and $t$ as above,
\[
\lambda^2 = k(s,t)=(s-t)\Bigg(\frac{1}{h(s)}-s+t\Bigg).
\]
Solving this equation for $t$, we get that $t$ is one of the two
numbers
\[
t_\pm = s - \frac{1}{h(s)}\pm \frac{\sqrt{1+4\lambda^2h(s)^2}}{2h(s)}.
\]
If $t=t_-$, then
\[
s-t>\frac{1}{2h(s)},
\]
and since $\frac{1}{h(s)}\rightarrow \infty$ as $s\rightarrow\infty$, this
can not hold for large $t$ because of (i). Hence,
$t=t_+$ for $t$ sufficiently large. $\endproof$

\vspace{.2cm}

Combining the previous lemmas we get:

\begin{prop}\label{prop5.13} Let $T\in\tilde\CM$ be an $R$--diagonal element, let
  $\lambda\in\C\setminus\{0\}$, and define $h(s)$ and $h_\lambda(s)$
  as in Lemma~\ref{lemma h}. Let $\mu= \mu_{|T|}$, and let
  $s(|\lambda|,t)$ be as in Definition~\ref{def5.10}. Then
  \[
  h_\lambda(s(|\lambda|,t))=h(t), \qquad t>0.
  \]
\end{prop}

\proof According to Lemma~\ref{lemma5.12}, if $t>t_{|\lambda|}$ and
$s=s(|\lambda|,t)$, then
\[
t=s+\frac{\sqrt{1+4\lambda^2h(s)^2}-1}{2h(s)}.
\]
Since $s(|\lambda|,t)>t$, we infer from Lemma~\ref{lemma h} that for
$t$ sufficiently large,
\[
h_\lambda(t)=h(s(|\lambda|,t)).
\]
Hence, by Lemma~\ref{lemma5.1} and Lemma~\ref{lemma5.8}, the same
formula holds for all $t>0$. $\endproof$

\vspace{.2cm}

\begin{lemma}\label{lemma5.14} Let $T$ be an unbounded $R$--diagonal element in $\CMD$,
  let $\lambda\in\C\setminus\{0\}$, and let $t>0$. With
  $\mu=\mu_{|T|}$ and $s(|\lambda|,t)$ as in Definition~\ref{def5.10}
  we then have:
  \begin{equation}\label{5.15}
    \Delta\big((T-\lambda\unit)\cc(T-\lambda\unit)+t^2\unit\big)=
    \frac{|\lambda|^2}{|\lambda|^2 + (s(|\lambda|,t)-t)^2} \Delta\big(T\cc
    T + s(|\lambda|,t)^2\unit\big).
  \end{equation}
\end{lemma}

\proof Since $T$ is $R$--diagonal, $T\underset{\ast\CD}{\sim} cT$ for
all $c\in\T$. Hence, the left--hand side of \eqref{5.15} depends only
on $|\lambda|$. It therefore suffices to consider only the case
$\lambda>0$. For $\lambda,t>0$, let
\[
H(t)=\textstyle{\frac12}\, \log\Delta(T\cc T+t^2\unit)
\]
and
\[
H_\lambda(t)= \textstyle{\frac12}\,\log\Delta\big((T-\lambda\unit)\cc
(T-\lambda\unit)+t^2\unit\big).
\]
Then with $\mu_\lambda = \mu_{|T-\lambda\unit|}$,
\[
H(t)= \frac12 \int_0^\infty \log(u^2+t^2)\,\d\mu(u),
\]
and
\[
H_\lambda(t)=\frac12 \int_0^\infty \log(u^2+t^2)\,\d\mu_\lambda(u).
\]
Since $T$ and $T-\lambda\unit$ belong to $\CMD$, $H$ and $H_\lambda$
take values in $\R$. Moreover, $H$ and $H_\lambda$ are differentiable
with derivatives $H'(t)=h(t)$ and $H_\lambda'(t)= h_\lambda(t)$. Also,
since $T\in \CMD$,
\begin{equation}\label{5.16}
  \lim_{t\rightarrow\infty}(H(t)-\log t) = \frac12
  \lim_{t\rightarrow\infty} \int_0^\infty \log\Bigg(1
  +\frac{u^2}{t^2}\Bigg)\,\d\mu(u) = 0,
\end{equation}
and similarly
\begin{equation}\label{5.17}
  \lim_{t\rightarrow\infty}(H_\lambda(t)-\log t)=0.
\end{equation}
Fix $\lambda>0$ and $t_0>0$. There is a constant $C$ such that
\[
H_\lambda(t)=\int_{t_0}^t h_\lambda(t')\,\d t' +C.
\]
Moreover, according to Proposition~\ref{prop5.13},
\[
h_\lambda(t)=h(s(\lambda,t)), \qquad t>0.
\]
Put $s(t)=s(\lambda,t)$ and $u(t)=t-s(t)$. Then $s(t)+u(t)=t$ and
$s'(t)+u'(t)=1$. Moreover, by Definition~\ref{def5.10},
\[
(s(t)-t)\Bigg(\frac{1}{h(s(t))}-s(t)+t\Bigg)=\lambda^2.
\]
Hence,
\[
u(t)\Bigg(\frac{1}{h(s(t))}-u(t)\Bigg)=\lambda^2,
\]
implying that
\[
h(s(t))=\frac{u(t)}{\lambda^2 + u(t)^2}.
\]
It follows that
\begin{eqnarray*}
  \int_{t_0}^t h_\lambda(v)\,\d v &=& \int_{t_0}^t
  h(s(v))(s'(v)+u'(v))\,\d v\\
  &=& \int_{t_0}^t \Bigg(h(s(v))s'(v) +\frac{u(v)}{\lambda^2 +
  u(v)^2}u'(v)\Bigg)\,\d v\\
&=& H(s(t))-H(s(t_0)) + \frac12 \log\Bigg(\frac{\lambda^2}{\lambda^2+u(t)^2}\Bigg)+\frac12\log\Bigg(\frac{\lambda^2+u(t_0)^2}{\lambda^2}\Bigg).
\end{eqnarray*}
Hence,
\[
H_\lambda(t)= H(s(t)) + \frac12 \log\Bigg(\frac{\lambda^2}{\lambda^2 +
  (s(t)-t)^2}\Bigg) + C',
\]
for a constant $C'$. 
Recall that $s(t)-t\rightarrow 0$ as $t\rightarrow \infty$
(cf. Lemma~\ref{lemma5.12}). It then follows from \eqref{5.16} and \eqref{5.17} that $C'$ must be
0. This finally shows us that
\[
\exp(2H_\lambda(t))= \frac{\lambda^2}{\lambda^2 +
  (s(t)-t)^2}\,\exp(2H(t)),
\]
and this proves \eqref{5.15}. $\endproof$

\vspace{.2cm}

\begin{thm}\label{thm5.15}
  Let $T\in\CMD$ be $R$--diagonal, let $\mu=\mu_{|T|}$, and let
  $s(|\lambda|,0)$ be as in Definition~\ref{def5.10}.
  \begin{itemize}
    \item[(i)] If $\lambda_1(\mu)<|\lambda|<\lambda_2(\mu)$, then
      \[
      \Delta(T-\lambda\unit)=
      \Bigg(\frac{|\lambda|^2}{|\lambda|^2+s(|\lambda|,0)^2}\,
      \Delta(T\cc T +s(|\lambda|,0)^2\unit)\Bigg)^\frac12.
      \]
    \item[(ii)] If $|\lambda|\leq \lambda_1(\mu)$, then
    $\Delta(T-\lambda\unit)= \Delta(T)$.
    \item[(iii)] If $|\lambda|\geq \lambda_2(\mu)$, then
    $\Delta(T-\lambda\unit)=|\lambda|$.
  \end{itemize}
\end{thm}

\proof The theorem is obviously true for $\lambda =0$. Moreover, as in
the proof of Lemma~\ref{lemma5.14}, it suffices to consider the case
$\lambda >0$. Note that
\begin{equation}\label{5.15a}
\Delta(T-\lambda\unit)^2 = \lim_{t\rightarrow 0+}
\Delta\big((T-\lambda\unit)\cc (T-\lambda\unit)+t^2\unit\big).
\end{equation}
Hence, (i) follows from Lemma~\ref{lemma5.8} and
Lemma~\ref{lemma5.14}. If $0<\lambda\leq \lambda_1(\mu)$, then by
Remark~\ref{rem5.11}, $\lim_{t\rightarrow 0+}s(\lambda,t)=0$. Hence,
(ii) also follows from Lemma~\ref{lemma5.14}. Now suppose $\lambda\geq
\lambda_2(\mu)$. Then $s(\lambda,t)\rightarrow \infty$ as
$t\rightarrow 0+$. The right--hand side of \eqref{5.15} is equal to
\[
\frac{\lambda^2 s(\lambda,t)^2}{\lambda^2 +
  (s(\lambda,t)-t)^2}\,\frac{\Delta(T\cc T
  -s(\lambda,t)^2\unit)}{s(\lambda,t)^2},
\]
where the first factor converges to $\lambda^2$ as $t\rightarrow 0+$,
and the second factor converges to 1 (cf. \eqref{5.16}). (iii) now follows
from \eqref{5.15} and \eqref{5.15a}. $\endproof$

\vspace{.2cm}

\begin{remark}\label{rem5.16} Note that
  \[
  \lambda_2(\mu)=\Bigg(\int_0^\infty u^2\,\d\mu_{|T|}(u)\Bigg)^\frac12
  = \|T\|_2
  \]
  and
  \[
  \lambda_1(\mu)= \Bigg(\int_0^\infty u^{-2}\,\d\mu_{|T|}(u)\Bigg)^{-\frac12}
  = \|T^{-1}\|_2^{-1},
  \]
  where $ \|T^{-1}\|_2:=+\infty$ in case ${\rm ker}(T)\neq 0$.
\end{remark}

\vspace{.2cm}

\begin{thm}\label{thm5.17}
  Let $T$ be an $R$--diagonal element in $\CMD$ with Brown measure
  $\mu_T$, and suppose $\mu_{|T|}$ is not a Dirac measure.
  \begin{itemize}
    \item[(a)] If ${\rm ker}(T)=0$, then
      \[
      {\rm supp}(\mu_T)=\{\lambda\in\C\,|\, \|T^{-1}\|_2^{-1}\leq
      |\lambda|\leq  \|T\|_2\}.
      \]
      Moreover, the $S$--transform of $\mu_{|T|^2}$ is well--defined
      and strictly increasing on $(-1,0)$ with
      \[
      \CS_{\mu_{|T|^2}}((-1,0))= (\|T\|_2^{-2},\|T^{-1}\|_2^2),
      \]
      and $\mu_T$ is the unique probability measure on $\C$ which
      is invariant under rotations and satisfies
      \[
      \mu_T\big(B(0, \CS_{\mu_{|T|^2}}(t-1)^{-\frac12})\big) = t,
      \quad 0<t<1.
      \]
    \item[(b)] If ${\rm ker}(T)\neq 0$, let $P$ denote the projection
    onto ${\rm ker}(T)$. Then
    \[
    {\rm supp}(\mu_T)=\{\lambda\in\C\,|\,|\lambda|\leq \|T\|_2\}.
    \]
    Moreover, the $S$--transform of $\mu_{|T|^2}$ is well--defined
    and strictly increasing on \mbox{$(\tau(P)-1,0)$} with
    \[
    \CS_{\mu_{|T|^2}}((\tau(P)-1,0))= (\|T\|_2^{-2},\infty),
    \]
    and $\mu_T$ is the unique probability measure on $\C$ which
      is invariant under rotations and satisfies
      \[
      \mu_T\big(B(0, \CS_{\mu_{|T|^2}}(t-1)^{-\frac12})\big) = t,
      \quad \tau(P)<t<1.
      \]
 \end{itemize}
\end{thm}

\proof By definition, $\d\mu_T(\lambda)=
\frac{1}{2\pi}\nabla^2\big(\log\Delta(T-\lambda\unit)\big)\d\lambda$ (in the
distribution sense). Hence, $\mu_T$ can be determined from
Theorem~\ref{thm5.15} in the same way as \cite[Theorem~4.4.]{HL} is
obtained from \cite[(4.5)]{HL}:

Using the same notation as in \cite{HL}, we define functions
$f,g:(0,\infty)\rightarrow \R$ by
\[
f(v)=\int_0^\infty \frac{1}{1+v^2w^2}\,\d\mu_{|T|}(w),
\]
and
\[
g(v)= \frac{1-f(v)}{v^2f(v)}.
\]
Moreover, for $\lambda\in(\|T^{-1}\|_2^{-2},\|T\|_2^2)$, let
$v(\lambda)$ denote the unique $v\in (0,\infty)$ such that
$g(v)=\lambda^2$. Then, in our notation,
\[
f(v)= \tau\big((\unit +v^2 T\cc T)^{-1}\big)= v^{-1}h(v^{-1}),
\]
and
\[
g(v)= v^{-1}\Bigg(\frac{1}{h(v^{-1})}-v^{-1}\Bigg)=k(v^{-1},0).
\]
Hence,
\[
v(\lambda)=\frac{1}{s(\lambda,0)},
\]
and it follows that the formula (4.15) in \cite{HL},
\[
\log\Delta(T-\lambda\unit)=\frac12 \int_0^\infty
\log(1+v^2w^2)\,\d\mu_{|T|}(w) +
\frac12\log\Bigg(\frac{\lambda^2}{1+v^2\lambda^2}\Bigg), \qquad
\lambda\in(\|T^{-1}\|_2^{-2},\|T\|_2^2),
\]
is equivalent to the one in Theorem~\ref{thm5.15}~(i). The rest of the
proof of Theorem~\ref{thm5.17} is identical to the second
part of the proof of \cite[Theorem~4.4]{HL}, since boundedness of $T$
is not a necessary assumption in the latter. $\endproof$

\vspace{.2cm}

\begin{remark} Let $T\in\CMD$ be $R$--diagonal. Then
  $\supp(\mu_T)\subseteq \sigma(T)$, and according to
  Theorem~\ref{thm5.17},
  \[
      {\rm supp}(\mu_T)=\{\lambda\in\C\,|\, \|T^{-1}\|_2^{-1}\leq
      |\lambda|\leq  \|T\|_2\}.
      \]
  Moreover, by arguments similar to the ones given in \cite[proof of
  Proposition~4.6]{HL}, one can show that
  \begin{itemize}
    \item[(a)] if $0<|\lambda|<\|T^{-1}\|_2^{-1}$, then
    $\lambda\in\sigma(T)$ iff $T$ does not have a bounded inverse, and
    \item[(b)] if $|\lambda|>\|T\|_2$, then $\lambda\in\sigma(T)$ iff
    $T$ is not bounded.
  \end{itemize}
\end{remark}
\section{Properties of $z=xy^{-1}$}

Let $\CM=L(\F_4)$ be the von Neumann algebra associated with the free
group on 4 generators. According to \cite{V1} or \cite{VDN}, $\CM$ is
a II$_1$--factor generated by a semicircular system
$(s_1,s_2,s_3,s_4)$, i.e. the $s_i$'s are freely independent
self--adjoint elements w.r.t. the unique tracial state $\tau$ on
$\CM$, and $s_i$ has distribution
\[
\d\mu_{s_i}(t)=\frac{1}{2\pi}\,\sqrt{4-t^2}\,1_{[-2,2]}(t)\,\d t,
\qquad 1\leq i\leq 4. 
\]
Put
\[
x=\frac{s_1+\i s_2}{\sqrt 2} \quad {\rm and} \quad y=\frac{s_3+\i s_4}{\sqrt 2}.
\]
Then $\CM=W\cc(x,y)$, and $(x,y)$ is a circular system in the sense of
\cite{VDN}. Also, by \cite{VDN}, $|y|$ has the distribution
\[
\d\mu_{|y|}(t)=\frac{2}{\pi}\,\sqrt{4-t^2}\,1_{[0,2]}(t)\,\d t.
\]
In particular, ${\rm ker}(y)=0$. In this section we will study the
unbounded operator
\[
z=xy^{-1}
\]
as well as its powers $z^n$, $n=2,3,\ldots$ We will need the following
simple observation:

\begin{lemma}\label{lemma6.1} Let $(\mu_n)_{n=1}^\infty$ and $\mu$ be
  probability measures on $\R$ with densities $(f_n)_{n=1}^\infty$ and
  $f$, respectively, w.r.t. Lebesgue measure. If $f_n\overset{n\rightarrow\infty}{\rightarrow}
  f$ a.e. w.r.t. Lebesgue measure, then $\mu_n\overset{n\rightarrow\infty}{\rightarrow} \mu$
  weakly.
\end{lemma}

\proof Recall that  $\mu_n\overset{n\rightarrow\infty}{\rightarrow} \mu$
  weakly iff for all $\phi\in C_0(\R)$,
  \begin{equation}\label{6.1}
    \lim_{n\rightarrow\infty}\int_\R \phi\,\d\mu_n =
    \int_\R\phi\,\d\mu.
  \end{equation}
Then, let $\phi\in C_0(\R)$ with $0\leq \phi\leq 1$. \eqref{6.1}
follows for such $\phi$ by application of Fatou's Lemma to each of the
sequences of integrals $\Bigg(\int_\R\phi f_n\d m\Bigg)_{n=1}^\infty$
and  $\Bigg(\int_\R(1-\phi)f_n\d m\Bigg)_{n=1}^\infty$. $\endproof$

\vspace{.2cm}

\begin{thm}\label{thm6.2} Let $(\CM,\tau)$ and $z=xy^{-1}$ be as
  above.
  \begin{itemize}
    \item[(a)] $z$ is an unbounded, $R$--diagonal operator.
    \item[(b)] The distribution of $z$ is given by
      \begin{equation}\label{6.2}
        \d\mu_{|z|}(t)=\frac2\pi
        \frac{1}{1+t^2}\,1_{(0,\infty)}(t)\,\d t.
      \end{equation}
    \item[(c)] For $p\in (0,1)$, $z, z^{-1}\in \Lp$, and
      \begin{equation}\label{6.3}
        \|z\|_p^p=\|z^{-1}\|_p^p= \Big[\cos\Big(\frac{p\pi}{2}\Big)\Big]^{-1}
        <\infty.
      \end{equation}
    \item[(d)] $z,z^{-1}\in\CMD$, and the Brown measure of $z$ is
    given by
    \begin{equation}\label{6.4}
      \d\mu_z(s)=\frac{1}{\pi(1+|s|^2)^2}\,\d s,
    \end{equation}
    where $\d s= \d \re s\,\d \im s$ is Lebesgue measure on $\C$.
  \end{itemize}
\end{thm}

\proof (a) Let $x=u|x|$ and $y=v|y|$ be the polar decompositions of $x$ and $y$. Then, according to \cite{VDN}, $u,|x|,v$ and $|y|$ are $\ast$--free elements, and $u$ and $v$ are Haar unitaries. In particular, $x$ and $y$ are $R$--diagonal and so is $y^{-1}$ (cf. Proposition~\ref{Rdiag1}). Moreover, $y^{-1}$ has polar decomposition
\[
y^{-1}=v\cc (v|y|^{-1}v\cc)=v\cc |y\cc|^{-1},
\]
which implies that $y^{-1}$ is affiliated with $W\cc(y)$. Hence, $x$ and $y^{-1}$ are $\ast$--free, and it follows from Proposition~\ref{Rdiag2} that $z=xy^{-1}$ is $R$--diagonal with
\[
\CS_{\mu_{|z|^2}}(t) = \CS_{\mu_{|x|^2}}(t)\,\CS_{\mu_{|y^{-1}|^2}}(t), \qquad t\in (-1,0).
\]
The distribution of $|x|^2$ has density
\[
\d\mu_{|x|^2}(t)=\frac{1}{2\pi}\sqrt{\frac{4-t}{t}}\,1_{(0,4)}(t)\,\d t,
\]
and thus $\CS_{\mu_{|x|^2}}$ is given by
\[
\CS_{\mu_{|x|^2}}(t)=\frac{1}{1+t}
\]
for all $t$ in a neighborhood of $(-1,0)$ (cf. \cite[example~5.2]{HL}). Since $|y^{-1}|=|y\cc|^{-1}\underset{\ast\CD}{\sim} |y|^{-1} \underset{\ast\CD}{\sim}|x|^{-1}$, we get from Proposition~\ref{S4} that
\[
\CS_{\mu_{|y^{-1}|^2}}(t)=\frac{1}{\CS_{\mu_{|x|^2}}(-1-t)}=-t, \qquad t\in(-1,0).
\]
Then
\begin{equation}\label{6.5}
\CS_{\mu_{|z|^2}}(t)=-\frac{t}{1+t}, \qquad t\in(-1,0),
\end{equation}
and
\[
\chi_{\mu_{|z|^2}}(t)=\frac{t}{1+t}\CS_{\mu_{|z|^2}}(t)=-\Bigg(\frac{t}{1+t}\Bigg)^2, \qquad t\in(-1,0).
\]
The inverse function of $\chi_{\mu_{|z|^2}}$ is then
\[
\psi_{\mu_{|z|^2}}(u)=\frac{-\sqrt{-u}}{1+\sqrt{-u}}, \qquad u\in(-\infty,0),
\]
and it follows that
\begin{equation}\label{6.6}
G_{\mu_{|z|^2}}(\lambda)=  \frac1\lambda\Big(1+\psi_{\mu_{|z|^2}}\Big(\textstyle{\frac1\lambda}\Big)\Big)= \frac{1}{\lambda-\sqrt{-\lambda}}, \qquad \lambda<0. 
\end{equation}
Let $\sqrt{w}$ denote the principal value of the square root of $w$ for $w\in\C\setminus(\infty,0]$. Then both sides of \eqref{6.6} are analytic in $\C\setminus[0,\infty)$. Thus, \eqref{6.6} holds for all $\lambda\in \C\setminus[0,\infty)$, and it follows that for $t>0$,
\begin{equation}\label{6.7}
-\frac1\pi\,\lim_{u\rightarrow 0+}\im G_{\mu_{|z|^2}}(t+\i u)=-\frac1\pi \,\im \Bigg(\frac{1}{t+\i\sqrt t}\Bigg) =\frac1\pi \frac{1}{\sqrt t(t+1)}.
\end{equation}
For $\beta\in(0,1)$,
\begin{equation}\label{6.8}
\int_0^\infty \frac{t^{\beta-1}}{1+t}\,\d t = \frac{\pi}{\sin(\beta\pi)},
\end{equation}
(cf. \cite[p.~592, formula~613]{Ha}). The right--hand side of \eqref{6.7} therefore defines the density of a probability measure, and then, by Lemma~\ref{lemma6.1}, the probability measures
\[
\frac1\pi\,\im G_{\mu_{|z|^2}}(t+\i u)\,\d t, \qquad u>0,
\]
converge weakly to 
\begin{equation}\label{6.9}
\frac1\pi \frac{1}{\sqrt t(t+1)}\,1_{(0,\infty)}(t)\,\d t,
\end{equation}
as $u\rightarrow 0+$. Hence, by the inverse Stieltjes transform, $\d\mu_{|z|^2}(t)$ is given by \eqref{6.9}, and then
\[
\d\mu_{|z|}(t)=\frac2\pi\frac{1}{1+t^2}\,1_{(0,\infty)}(t)\,\d t.
\]
This proves (a) and (b). 

In order to prove (c), note that according to \eqref{6.8},
\begin{eqnarray*}
\tau(|z|^p)&=&\frac2\pi\int_0^\infty\frac{t^p}{1+t^2}\,\d t\\
&=& \frac1\pi \int_0^\infty \frac{w^{\frac{p-1}{2}}}{1+w}\,\d w\\
&=& \Big[\sin\Big(\frac{\pi(p+1)}{2}\Big)\Big]^{-1},
\end{eqnarray*}
proving (c). Since $\Lp\subseteq \CMD$, $p>0$, $z,z^{-1}\in\CMD$. According to Theorem~\ref{thm5.17}, $\mu_z$ is then the unique probability measure on $\C$ which is invariant under rotations and satisfies
\[
\mu_z\big(B(0,\CS_{\mu_{|z|^2}}(t-1)^{-\frac12})\big) = t, \qquad 0<t<1.
\]
Then by \eqref{6.5},
\[
\mu_z\Big(B\Big(0,\sqrt{\textstyle{\frac{t}{1-t}}}\Big)\Big) = t, \qquad 0<t<1,
\]
that is,
\[
\mu_z(B(0,r))=\frac{r^2}{1+r^2},\qquad r>0.
\]
Hence, $\frac{\d}{\d r}\mu_z(B(0,r))=\frac{2r}{(1+r^2)^2}$, and combining this with the fact that $\mu_z$ is invariant under rotations, we find that $\mu_z$ has density w.r.t. Lebesgue measure on $\C$ given by
\[
\frac{1}{2\pi r} \frac{2r}{(1+r^2)^2}= \frac1\pi \frac{1}{(1+r^2)^2}, \qquad r>0,
\]
where $r=|s|$, $s\in\C\setminus\{0\}$. This proves (d). $\endproof$

\vspace{.2cm}

\begin{lemma}\label{lemma6.3} Let $\mu$ be a probability measure on $[0,\infty)$ and, as in section~5, put
\[
h(s)=\int_0^\infty \frac{s}{s^2+u^2}\,\d\mu(u),\qquad s\in (0,\infty).
\]
Then for $0<p<2$,
\begin{equation}\label{6.10}
\int_0^\infty u^{-p}\,\d\mu(u)=\frac2\pi \sin\Big(\frac{\pi p}{2}\Big)\int_0^\infty s^{-p}h(s)\,\d s.
\end{equation}
\end{lemma}

\proof By Tonelli's theorem,
\[
\int_0^\infty s^{-p}h(s)\,\d s = \int_0^\infty\Bigg(\int_0^\infty\frac{s^{1-p}}{s^2+u^2}\,\d s\Bigg)\d\mu(u).
\]
Letting $s=ut^{\frac12}$, we find (using \eqref{6.8}) that
\[
\int_0^\infty\frac{s^{1-p}}{s^2+u^2}\,\d s = \frac12 u^{-p}\int_0^\infty\frac{t^{-\frac p2}}{1+t}\,\d t =\frac\pi2\Big[\sin\Big(\frac{\pi p}{2}\Big)\Big]^{-1}u^{-p}.
\]
This proves \eqref{6.10}. $\endproof$

\vspace{.2cm}

\begin{thm}\label{thm6.4} Let $(\CM,\tau)$ and $z$ be as in Theorem~\ref{thm6.2}, and let $n\in\N$. 
\begin{itemize}
\item[(a)] $z^n$ is an unbounded $R$--diagonal operator.
\item[(b)] \begin{equation}\label{6.11}
\int_0^\infty\frac{s}{s^2+u^2}\,\d\mu_{|z|^n}(u)=\Big(s+s^{\frac{n-1}{n+1}}\Big)^{-1}, \qquad s>0.
\end{equation}
\item[(c)] For $p\in \Big(0,\frac{2}{n+1}\Big)$, $z^n$ and $z^{-n}$ both belong to $\Lp$, and
\begin{equation}\label{6.12}
\|z^n\|_p^p=\|z^{-n}\|_p^p=\frac{(n+1)\sin\big(\textstyle{\frac{\pi p}{2}}\big)}{\sin\big(\textstyle{\frac{(n+1)\pi p}{2}}\big)}.
\end{equation}
\item[(d)] If $p\in \Big(0,\frac{2}{n+1}\Big)$ and $\lambda\in\C$, then ${\rm ker}(z^n-\lambda\unit)\neq 0$. Moreover, $(z^n-\lambda\unit)^{-1}\in\Lp$ with
\begin{equation}\label{6.13}
\|(z^n-\lambda\unit)^{-1}\|_p\leq \|z^{-n}\|_p.
\end{equation}
\end{itemize}
\end{thm}

\proof According to Proposition~\ref{Rdiag4}, $z^n$ is $R$--diagonal. Moreover, since
\[
\CS_{\mu_{|z|^2}}(t)^n=\Bigg(-\frac{t}{1+t}\Bigg)^n, \qquad t\in(-1,0),
\]
\[
\chi_{\mu_{|z^n|^2}}(t)=\frac{1}{1+t}\CS_{\mu_{|z^n|^2}}(t)=-\Bigg(-\frac{t}{1+t}\Bigg)^{n+1}, \qquad t\in(-1,0),
\]
with inverse function
\[
\psi_{\mu_{|z^n|^2}}(u)=-\frac{(-u)^{\frac{1}{n+1}}}{1+(-u)^{\frac{1}{n+1}}}, \qquad u\in(-\infty,0).
\]
Hence, for $\lambda\in (-\infty,0)$,
\begin{equation}\label{6.14}
G_{\mu_{|z^n|^2}}(\lambda)=\frac1\lambda\Big(1+\psi_{\mu_{|z^n|^2}}\big(\textstyle{\frac1\lambda}\big)\Big)=\frac{1}{\lambda\Big(1+(-\lambda)^{-\frac{1}{n+1}}\Big)}.
\end{equation}
Let
\[
h_n(s)=\int_0^\infty \frac{s}{s^2+u^2}\,\d\mu_{|z^n|}(u),\qquad s\in(0,\infty).
\]
Then
\begin{eqnarray*}
h_n(s)&=& s\,\tau\big((s^2\unit + |z^n|^2)^{-1}\big)\\
&=& -s\,G_{\mu_{|z^n|^2}}(-s^2)\\
&\overset{\eqref{6.14}}{=}& \Big(s+s^{\frac{n-1}{n+1}}\Big)^{-1}.
\end{eqnarray*}
This proves (b). 

Since $z=xy^{-1}$, where $(x,y)$ is a circular family, it is clear that $z^{-n}\underset{\ast\CD}{\sim}z^n$ for all $n\in\N$. Hence, $\|z^n\|_p=\|z^{-n}\|_p$ for all $p>0$. Note that for $p>0$,
\[
\|z^{-n}\|_p^p = \tau(|z^{-n}|^p)=\tau(|(z^n)\cc|^{-p})=\tau(|z^n|^{-p}).
\]
Thus, by Lemma~\ref{lemma6.3}, for $p\in(0,2)$,
\begin{equation}\label{6.15}
\|z^{-n}\|_p^p=\int_0^\infty u^{-p}\,\d\mu_{|z^n|}(u)=\frac2\pi\sin\Big(\frac{\pi p}{2}\Big)\int_0^\infty s^{-p}h_n(s)\,\d s.
\end{equation}
By application of \eqref{6.11} we find that
\[
\int_0^\infty s^{-p}h_n(s)\,\d s = \int_0^\infty \frac{s^{-p-\frac{n-1}{n+1}}}{s^{\frac{2}{n+1}}+1}\,\d s= \frac{n+1}{2}\int_0^\infty \frac{t^{-\frac{(n+1)p}{2}}}{1+t}\,\d t.
\]
Then by \eqref{6.15} and \eqref{6.8}, for $0<p<\frac{2}{n+1}$,
\begin{equation}
\begin{split}
\|z^{-n}\|_p^p &= (n+1)\sin\Big(\frac{\pi p}{2}\Big)\Big[\sin\Big(\pi\Big(1-\frac{(n+1)p}{2}\Big)\Big)\Big]^{-1}\\
& = (n+1)\sin\Big(\frac{\pi p}{2}\Big)\Big[\sin\Big(\frac{(n+1)\pi p}{2}\Big)\Big]^{-1},\label{6.16}
\end{split}
\end{equation}
and this proves (c). Note that the right--hand side of \eqref{6.16}
converges to $\infty$ as $p\rightarrow \frac{2}{n+1}-$. Hence,
$z^{-n}\notin L^{\frac{2}{n+1}}(\CM,\tau)$, and the same holds for
$z^n$. In particular, $z^n$ is not bounded, and this proves (a). In
order to prove (d), let $\lambda\in\C\setminus\{0\}$, and put
\[
h_{n,\lambda}(t)=\int_0^\infty\frac{t}{t^2+u^2}\,\d\mu_{|z^n-\lambda\unit|}(u), \qquad t>0.
\]
Then by Proposition~\ref{prop5.13},
\[
h_{n,\lambda}(t)=h_n(s_n(|\lambda|,t)), \qquad t>0,
\]
where $s_n(|\lambda|,t)$ is given by Definition~\ref{def5.10} in the
case $\mu=\mu_{|z^n|}$. Note that, according to
Definition~\ref{def5.10},
\[
s_n(|\lambda|,t)>t, \qquad t>0.
\]
Moreover, by \eqref{6.11}, $h_n$ is monotonically decreasing on
$(0,\infty)$. Thus,
\[
h_{n,\lambda}(t)\leq h_n(t), \qquad t>0.
\]
It now follows from Lemma~\ref{lemma6.3} that for $p\in(0,2)$,
\begin{equation}\label{6.17}
\int_0^\infty u^{-p}\,\d\mu_{|z^n-\lambda\unit|}(u)\leq \int_0^\infty u^{-p}\,\d\mu_{|z^n|}(u).
\end{equation}
According to (c), the right--hand side of \eqref{6.17} is finite for
$p\in \Big(0,\frac{2}{n+1}\Big)$. Hence, for such $p$, ${\rm
  ker}(z^n-\lambda\unit)=0$, $(z^n-\lambda\unit)^{-1}\in\Lp$, and
\[
\|(z^n-\lambda\unit)^{-1}\|_p^p\leq \|z^{-n}\|_p^p. \qquad \endproof
\]

\vspace{.2cm}

\begin{remark}
Note that Theorem~\ref{thm6.4} (a) and (c) generalize Theorem~\ref{thm6.2} (a) and
(c) to  all $n\in\N$. It is not hard to generalize
Theorem~\ref{thm6.2} (b) and (d) as well. One finds that the
distribution of $|z^n|$ is given by
\[
\d\mu_{|z^n|}(t)=\frac2\pi
\frac{\sin\Big(\frac{\pi}{n+1}\Big)}{t\Big(t^{\frac{2}{n+1}} +
  2\cos\Big(\frac{\pi}{n+1}\Big)+t^{-\frac{2}{n+1}} \Big)}\,1_{(0,\infty)(t)}\,\d t,
\]
and the Brown measure of $z^n$ is given by
\[
\d\mu_{z^n}(s)= \frac{1}{n\pi}\frac{|s|^{\frac2n
    -2}}{(1+|s|^{\frac2n})^2}\,\d\re s\,\d\im s.
\]
We leave the details of proof to the reader. 
\end{remark}

{\small

\noindent Uffe~Haagerup\\
Department of Mathematics and Computer Science\\
University of Southern Denmark\\
Campusvej~55, 5230~Odense~M\\
Denmark\\
{\tt haagerup@imada.sdu.dk}

\vspace{.5cm}

\noindent Hanne~Schultz\\
Department of Mathematics and Computer Science\\
University of Southern Denmark\\
Campusvej~55, 5230~Odense~M\\
Denmark\\
{\tt schultz@imada.sdu.dk}

\end{document}